\documentclass[11pt]{article}

\usepackage[english]{babel}
\usepackage[T1]{fontenc}
\usepackage[utf8]{inputenc}

\usepackage[top=2cm, bottom=2cm, outer=3cm, inner=2cm,
            heightrounded, marginparwidth=2.5cm, marginparsep=1cm]{geometry}

\usepackage{amsmath, amsthm, amssymb}
\usepackage{enumerate}
\usepackage{tikz}
\usetikzlibrary{decorations.pathmorphing}

\usepackage{caption}
\usepackage{float}
\usepackage{framed}
\usepackage{hyperref}
\usepackage{xspace}

\newtheorem{theorem}{Theorem}[section]
\newtheorem{proposition}[theorem]{Proposition}
\newtheorem{lemma}[theorem]{Lemma}

\newtheorem{corollary}[theorem]{Corollary}
\newtheorem{observation}[theorem]{Observation}

\newtheorem{definition}[theorem]{Definition}

\newcommand{\secref}[1]{Section~\ref{#1}\xspace}
\newcommand{\thmref}[1]{Theorem~\ref{#1}\xspace}
\newcommand{\lemref}[1]{Lemma~\ref{#1}\xspace}
\newcommand{\propref}[1]{Proposition~\ref{#1}\xspace}
\newcommand{\corref}[1]{Corollary~\ref{#1}\xspace}

\newcommand{\obsref}[1]{Observation~\ref{#1}\xspace}

\newcommand{\figref}[1]{Figure~\ref{#1}\xspace}
\newcommand{\rt}{\right}
\newcommand{\lt}{\left}

\newif\ifnotescoauthor
\notescoauthortrue 
\newcommand{\comment}[2]{\ifnotescoauthor \begin{framed} $\blacktriangleright$
{\sf #1} $\blacktriangleleft$ \ifx&#2&\else\medskip\hrule\medskip
#2\fi\end{framed}\fi}

\newcommand{\block}[1]{#1_{\textup{block}}}
\newcommand{\Bin}{\textup{Bin}}
\newcommand{\aas}{a.a.s.\ }
\newcommand{\gnp}{\mathbb{G}}
\newcommand{\Step}[1][Step]{\underline{\emph{#1:}}\enskip}
\newcommand{\Prob}{\mathbb{P}}
\newcommand{\E}{\mathbb{E}}
\newcommand{\eps}{\varepsilon}
\newcommand{\mc}[1]{\mathcal #1}

\AtBeginDocument{\setlength{\abovedisplayskip}{5pt}}
\AtBeginDocument{\setlength{\belowdisplayskip}{5pt}}
\newcommand{\smalllist}{
\setlength{\itemsep}{1pt}
\setlength{\parskip}{0pt}
\setlength{\parsep}{0pt}}

\title{The threshold probability for long cycles}
\author{
Roman Glebov\thanks{Department of Mathematics, ETH, 8092 Zurich, Switzerland.
Email: \texttt{roman.glebov@math.ethz.ch}.} \and
Humberto Naves\thanks{Department of Mathematics, ETH, 8092 Zurich, Switzerland
and Department of Mathematics, UCLA, Los Angeles, CA 90095 USA.
Email: \texttt{hnaves@math.ucla.edu}.} \and
Benny Sudakov\thanks{Department of Mathematics, ETH, 8092 Zurich, Switzerland.
Email: \texttt{benjamin.sudakov@math.ethz.ch}. Research supported in part by
SNSF grant 200021-149111 and by a USA-Israel BSF
grant.}}

\date{}

\begin{document}
\maketitle
\setcounter{page}{1}

\begin{abstract}

  For a given graph $G$ of minimum degree at least $k$, let
  $G_p$ denote the random spanning subgraph of $G$ obtained by
  retaining each edge independently with probability $p=p(k)$. We prove
  that if $p \ge \frac{\log k + \log \log k + \omega_k(1)}{k}$,
  where $\omega_k(1)$ is any function tending to infinity with $k$,
  then $G_p$ asymptotically almost surely contains a cycle
  of length at least $k+1$. When we take $G$ to be the complete
  graph on $k+1$ vertices, our theorem coincides with the classic
  result on the threshold probability for the existence of a
  Hamilton cycle in the binomial random graph.

\end{abstract}

\section{Introduction}

Given a graph $G$ and a real $p\in [0, 1]$, let $G_p$ be the probability
space of subgraphs of $G$ obtained by taking each edge of $G$ independently
with probability $p$. We sometimes use the notation $(G)_p$ to
avoid ambiguity. For a given graph property $\mc P$ and sequences
of graphs $\{G_i\}_{i=1}^\infty$ and probabilities $\{p_i\}_{i=1}^\infty$,
we say that $(G_i)_{p_i}\in \mc P$ \emph{asymptotically almost surely},
or \aas for brevity, if the probability that $(G_i)_{p_i}\in\mc P$
tends to $1$ as $i$ goes to infinity. In this paper, when $G$ and $p$
depend upon some parameter, we abuse notation and consider $G$ and $p$
as sequences obtained by taking the parameter to tend to infinity, and we
say that $G_p$ has $\mc P$ \aas if the sequence does.

When the host graph $G$ is the complete graph on $n$ vertices, the random
graph model $G_p$ coincides with the classic binomial random graph
model $\gnp(n,p)$, introduced independently by Gilbert in \cite{Gilbert}
and by Erd\H{o}s and R\'{e}nyi in \cite{ErdRen}. This important model
has been studied extensively for the past few decades. A result of P\'{o}sa~%
\cite{Posa} states that for some large constant $C>0$, if $p \ge
\frac{C\log n}{n}$  then $\gnp(n,p)$ \aas contains a  Hamilton cycle.
This result was later strengthened by Korshunov~\cite{Korshunov}, Koml\'{o}s
and Szemer\'{e}di~\cite{KomSze}, and independently by
Bollob\'{a}s~\cite{Bollobas}. They proved that the same statement holds
for $p \ge\frac{\log n + \log \log n + \omega_n(1)}{n}$, provided $n$
is large.

In this paper we extend the aforementioned result to a more general
class of graphs. More precisely, we would like to replace the host graph $G$,
taken to be the complete graph in the classic setting, by a graph with
minimum degree at least $k$, and to find \aas a cycle of length at least $k+1$
in the random subgraph $G_p$. Our main result is as follows.

\begin{theorem}

  \label{thm:main}
  Let $G$ be a graph with minimum degree at least $k$. If
  $p=p(k)\ge\frac{\log k + \log \log k + \omega_k(1)}{k}$,
  then $G_p$ \aas contains a cycle of length at least $k+1$.

\end{theorem}

Our results are complimentary to the ones of
Krivelevich, Lee, and Sudakov~\cite{KrivLeeSud} and of Riordan~\cite{Riordan}.
They proved that for $p = \frac{\omega_k(1)}{k}$, the graph $G_p$ \aas contains a
cycle of length at least $(1+o(1))k$, which might be slightly less than
$k+1$. Since the property stated in the main theorem is monotone increasing,
we may assume throughout the paper that $p\le\frac{\log k+2\log\log k}{k}$.

The rest of this paper is organized as follows. \secref{sec:preliminaries}
contains a variety of tools, which are used to prove \thmref{thm:main}.
All propositions, statements and lemmas in that section are stated without
proofs. In \secref{sec:proof}, we prove our main theorem. The final section
contains some concluding remarks.

\subsection{Notation}

A graph $G=(V,E)$ is given by a pair of its (finite) vertex set $V(G)$ and
edge set $E(G)$. We use $|G|$ or $|V(G)|$ to denote the order of the
graph. For a subset $X$ of vertices, we use $e(X)$ to denote the number
of edges spanned by $X$, and for two disjoint sets $X,Y$, we use $e(X,Y)$ to
denote the number of edges with one endpoint in $X$ and the other in $Y$.
Let $G[X]$ denote the subgraph of $G$ induced by a subset of
vertices $X$. We write $N(X)$ to denote the collection of vertices
outside of $X$ that have at least one neighbor in $X$.
When $X$ consists of a single vertex, we abbreviate $N(v)$
for $N(\{v\})$, and let $\deg(v)$ denote the cardinality of $N(v)$, i.e.,
the degree of $v$. For two graphs $G_1$ and $G_2$, not necessarily over the
same vertex set, we define their intersection as $G_1\cap G_2 =
(V(G_1) \cap V(G_2), E(G_1)\cap E(G_2))$, and union as $G_1\cup G_2 =
(V(G_1) \cup V(G_2), E(G_1)\cup E(G_2))$.
Moreover, if $X$ is a set of vertices, we let $G\setminus X$ to be the
induced subgraph $G[V(G)\setminus X]$. Finally, if $G$ is a graph
and $E$ is a collection of unordered pairs of vertices from $V(G)$,
let $G+E$ denote the graph obtained from $G$ by adding the edges in $E$
which are not already in $G$. When there are several graphs
under consideration, we use subscripts such as $N_G(X)$ indicating
the relevant graph of interest.

The probability space $G_p$ is a simple product space.
When sampling from this model, one could unveil the graph $G_p$ by
successively answering queries of the form ``does $e$ belong to $G_p$?''
for each edge $e\in E(G)$. Since the answers to these queries are
independent, this process can be carried out regardless of the order
of the queries, as long as each edge of $G$ is queried exactly once.
Throughout the paper we expose $G_p$ in this manner.
The edges of $G$ not yet queried in $G_p$ shall be named \emph{untested},
while the others are called \emph{tested}. When an edge $e$ from $G$
is queried and the outcome turns out to be positive, we say that $e$ was
\emph{successfully tested}, or equivalently, $e$ was
\emph{successfully exposed}. We write \emph{partially exposed} $G_p$
as a reminder that not all edges of $G$ were tested in $G_p$. All
probabilistic statements involving a partially exposed $G_p$ must
be conditioned on the outcome of the tested edges at that particular
moment of the exposure process.
More precisely, if $Q$ is the set of testes edges, and $E\subseteq Q$
is the set of successfully tested edges of the partially exposed
$G_p$, then for each subgraph $\Gamma\subseteq G$, the probability
that we obtain the graph $\Gamma$ after we expose all the remaining
untested edges is $\Prob\big[G_p = \Gamma ~\Big|~ E(G_p) \cap Q = E\big]$.

To simplify the presentation, we often omit
floor and ceiling signs whenever these are not crucial and make no
attempts to optimize the absolute constants involved. We also assume
that the parameter $k$ (which always denotes the minimum degree of the
host graph) tends to infinity and therefore is sufficiently
large whenever necessary. All our asymptotic notation symbols
($O$, $o$, $\Omega$, $\omega$, $\Theta$) are relative to this
variable $k$, unless otherwise specified with a subscript.
Finally, all logarithms are to base $e\approx 2.718$.

\section{Preliminaries}
\label{sec:preliminaries}

\subsection{Probabilistic tools}
\label{sec:prob_tools}

We use extensively the following well-known bounds on the lower and
upper tails of the binomial distribution due to Chernoff (see,
e.g., \cite[Theorems A.1.11, A.1.13, and A.1.12]{AlonSpencer}).

\begin{lemma}

  \label{lem:chernoff}
  If $X \sim \Bin(n,p)$, then
  \begin{itemize} \smalllist
    \item $\Prob\left[X<(1-a)np\right]<\exp\left(-\frac{a^2np}{2}\right)$
          for every $a>0.$
    \item $\Prob\left[X>(1+a)np\right]<\exp\left(-\frac{a^2np}{3}\right)$
          for every $0 < a < 1.$
  \end{itemize}

\end{lemma}

\begin{lemma}

  \label{lem:chernoff2}
  Let $X \sim \Bin(n,p)$ and $a \in \mathbb{N}$. Then
  $\Prob[X \ge a] \leq \left(\frac{enp}{a}\right)^a$.

\end{lemma}

\subsection{Depth-First Search algorithm}
\label{sec:dfs}
Depth-First Search (DFS) is a well-known graph exploration algorithm,
usually applied to discover the connected components of an input graph.
The algorithm visits all vertices of a graph $H$ (the input of the DFS)
and produces a rooted spanning forest $T$ of $H$ (the output). It also
maintains a stack $S$ (last-in-first-out data structure) of vertices.
Initially, the stack is empty, and all vertices of $H$ are \emph{active}.
Each active vertex $v$ eventually gets reached, henceforth becoming
inactive, and is then pushed into $S$. At some point later, the same
vertex $v$ is popped from $S$ and is declared \emph{explored}.
Once a vertex becomes explored, it never changes its state back to active
again. Indeed, the algorithm ends when all the vertices of $H$ become explored.
The main loop of the DFS is as follows.
\begin{enumerate}[(i)] \smalllist
\item If $S$ is empty, choose an active vertex $v$, deactivate it, and push it
onto the stack. The vertex $v$ is the root of a new tree in $T$.
\item Otherwise, let $u$ be the unique vertex on top of the stack $S$. The
algorithm then queries for active neighbors of $u$ in $H$, i.e., active
vertices $w$ such that $uw$ forms an edge in $H$. If there is such an edge, we
remove $w$ from the set of active vertices and place it on top of $S$.
Otherwise, we just pop $u$ from the top of $S$ and mark it as explored.
\end{enumerate}
Notice that we specified neither how to choose the new vertex $v$ in (i) nor
the order in which the neighbors of $u$ should be queried in (ii).
It was implicitly assumed that these choices were made
according to some predetermined order --- the \emph{priority} of the
DFS.

The rooted spanning forest $T$ produced by the DFS induces a partial
order on the vertices of $H$. Namely, we say that $u \le_T v$ if $u$
belongs to the (unique) path connecting $v$ to a root of $T$. In this
case, we say that $u$ is an \emph{ancestor} of $v$, or equivalently,
$v$ is a \emph{descendant} of $u$ with respect to $T$. Whenever
$uv\in E(T)$, we say that $v$ is an \emph{immediate descendant} of $u$,
or, equivalently, $u$ is an \emph{immediate ancestor} of $v$.
A key observation is the following.

\begin{proposition}

  \label{prop:dfs_tree}
  For every edge $uv$ of $H$, $u$ and $v$ are comparable with respect to
  $\le_T$.

\end{proposition}

In our setting, we utilize the DFS algorithm on the random graph $G_p$,
and expose an edge only at the moment when its existence is queried by
the algorithm. Note that the input graph $G_p$ might be already
partially exposed at the moment we start the DFS. In this
case it is perfectly possible that the algorithm reuses some of
the successfully exposed edges (the algorithm never queries
the same edge twice). We discuss this topic in more detail in
\secref{sec:step2}.

Regardless of the portion of $G_p$ that was already exposed, the
following is always true.
\begin{proposition}

  \label{prop:dfs_new_edge}
  The rooted forest $T$ produced by the DFS algorithm running on
  a partially exposed $G_p$ contains all successfully tested edges
  revealed by the algorithm.

\end{proposition}
For instance, if we apply the DFS to $G_p$ with all the edges of $G$
initially untested, since the resulting forest $T$ has at most $n-1$ edges,
the algorithm must necessarily stop after the first $n-1$ successfully
exposed edges. Moreover, the connected components of $G_p$, when viewed as
vertex subsets of $V(G)$, are the same as the components of $T$, regardless
of the outcome of the remaining untested edges from $G_p$. One
noteworthy advantage of the DFS algorithm is that it produces this
``certificate'' for the connected components of a random graph by testing
very few of its edges. For more details on the application of the
depth-first search algorithm to random graphs, we refer the reader to
\cite{KrivSud}.

\subsection{Block algorithm}
\label{sec:block_algorithm}

Let us briefly recall some standard definitions and notions in
graph theory. Let $H$ be a graph. A vertex in $H$ is a \emph{cut-vertex}
if by removing it, we increase the number of connected components of $H$.
A maximal connected subgraph of $H$ without a cut-vertex is called a
\emph{block}. A \emph{$2$-connected graph} is a graph of order at least $3$
having no cut-vertex. In general, a \emph{$t$-connected graph} is a graph
$H$ of order at least $t+1$ such that $H\setminus X$ is connected for all
subsets $X\subseteq V(H)$ of size smaller than $t$.

In our quest to find long cycles, we will need to merge some already
revealed cycles into longer ones. To merge two disjoint cycles,
we need to find a collection of vertex-disjoint paths connecting them.
A classic result of Menger~\cite{Menger} enables us to find these paths.

\begin{theorem}[Menger]

  \label{thm:menger}
  Let $H$ be a $t$-connected graph. For every pair of subsets $A$ and
  $B$ of $V(H)$, there are at least $min\{t, |A|, |B|\}$ vertex-disjoint
  paths in $H$ that connect $A$ and $B$.

\end{theorem}

We extensively apply Menger's result inside the blocks of $G_p$.
This can be done because a block having at least $3$ vertices is necessarily
$2$-connected. To discover the blocks of $G_p$ we use another algorithm.
Our proposed algorithm produces a similar ``certificate'' for the
blocks of $G_p$, just like the DFS does for the connected components of $G_p$.

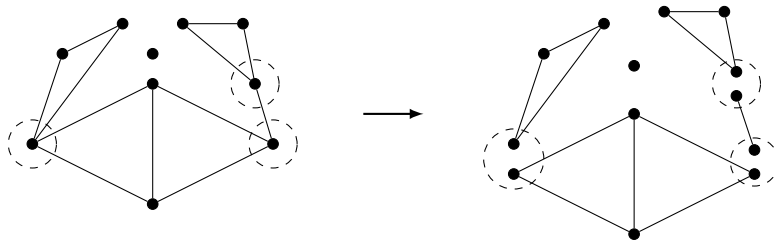
\begin{figure}[ht!]
\centering
\begin{tikzpicture}
  [scale=0.8, auto=left, every node/.style={circle, draw, fill=black,
   inner sep=0pt, minimum width=4pt}]
  \node (n1) at ( -1,  1.0 cm) {};
  \node (n2) at ( -1, -1.0 cm) {};
  \node (n3) at (-3,  0.0 cm) {};
  \node (n4) at ( 1,  0.0 cm) {};
  \node (n5) at (-1.5,  2.0 cm) {};
  \node (n6) at (-2.5,  1.5 cm) {};

  \node (n7) at (-0.5,  2.0 cm) {};
  \node (n8) at (0.5,  2.0 cm) {};
  \node (n9) at (0.7,  1.0 cm) {};
  \node (n10) at (-1.0, 1.5 cm) {};

  \node (m1) at ( 7,  0.5 cm) {};
  \node (m2) at ( 7, -1.5 cm) {};
  \node (m31) at (5,  -0.5 cm) {};
  \node (m32) at (5,  0.0 cm) {};
  \node (m41) at ( 9,  -0.5 cm) {};
  \node (m42) at ( 9,  -0.1 cm) {};
  \node (m5) at (6.5,  2.0 cm) {};
  \node (m6) at (5.5,  1.5 cm) {};

  \node (m7) at (7.5,  2.2 cm) {};
  \node (m8) at (8.5,  2.2 cm) {};
  \node (m91) at (8.7,  0.8 cm) {};
  \node (m92) at (8.7,  1.2 cm) {};
  \node (m10) at (7.0, 1.3 cm) {};

  \foreach \from/\to in {n1/n3, n1/n4, n2/n3, n2/n4, n1/n2,
                         n3/n5, n3/n6, n5/n6, n4/n9, n9/n7,
                         n9/n8, n7/n8}
    \draw (\from) -- (\to);

  \foreach \from/\to in {m1/m31, m1/m41, m2/m31, m2/m41, m1/m2,
                         m32/m5, m32/m6, m5/m6, m42/m91, m92/m7,
                         m92/m8, m7/m8}
    \draw (\from) -- (\to);

  \draw[style=dashed] (-3,0) circle (0.4);
  \draw[style=dashed] ( 1,0) circle (0.4);
  \draw[style=dashed] (0.7,1.0) circle (0.4);

  \draw[style=dashed] (5,-0.25) circle (0.5);
  \draw[style=dashed] (9,-0.3) circle (0.4);
  \draw[style=dashed] (8.7,1.0) circle (0.4);

  \draw[-latex, thick, ] (2.5, 0.5) -- (3.5, 0.5);
\end{tikzpicture}
\caption{The block decomposition of a graph.}
\label{fig:block_decomposition}
\end{figure}

Unlike the connected components of a graph $H$, the blocks of $H$ must not
necessarily be disjoint, as \figref{fig:block_decomposition} shows.
In fact two blocks can intersect, but in at most one vertex.
Moreover, it is well-known that blocks form a forest-like structure.
More formally, let $\block{H}$ be the bipartite graph on the vertex
set $\mc A\cup \mc B$, where $\mc A$ is the set of all
cut-vertices of $H$, $\mc B$ is the set of all blocks of $H$,
and the edges are formed by pairs $\{v,B\}$ satisfying $v\in \mc A$,
$B\in \mc B$ and $v \in B$. The resulting graph $\block{H}$,
referred to as the \emph{block decomposition of $H$}, is always cycle free.
This graph is also commonly known as the \emph{block-cutpoint graph of $H$}.

We summarize some of the properties of the block decomposition in the
next proposition. For more details, we refer the interested reader to
\cite[Chapter 3]{Diestel} and \cite[Chapter 4]{West}.

\begin{proposition}

  \label{prop:block_decomposition}
  Let $H$ be a graph and let $\block{H}$ be its block decomposition
  with vertex set $\mc A \cup \mc B$.

  \begin{enumerate}[(i)] \smalllist
  \item The equality $\bigcup_{B\in \mc B} V(B) = V(H)$ holds,
  and for every two distinct blocks $B,B'\in\mc B$, their intersection
  $B\cap B'$ is either empty or contains exactly one cut-vertex from
  $\mc A$. Furthermore, we have $|\mc B| \le |V(H)|$.

  \item The sets $E(B)$ for $B\in \mc B$ form a partition of $E(H)$.

  \item The graph $\block{H}$ is always cycle free. Moreover, $\block{H}$
  is a tree if $H$ is connected.
  \newcounter{tempCounter}
  \setcounter{tempCounter}{\value{enumi}}
  \end{enumerate}
  Furthermore, if $H$ and $H^*$ are two graphs having the same number
  of connected components, where $H$ spanning subgraph of $H^*$,
  then the following statements hold.
  \begin{enumerate}[(i)] \smalllist  \setcounter{enumi}{\value{tempCounter}}
  \item Every cut-vertex from $H^*$ is also a cut-vertex in $H$.
  \item If $v$ is a cut-vertex from $H$ but not a cut-vertex
  from $H^*$ then there exists an edge $e\in E(H^*)\setminus E(H)$
  which is not contained in any block of $H$.
  \item If $H$ and $H^*$ have the same set of cut-vertices then
  $\block{H} \simeq \block{H^*}$.
  \end{enumerate}

\end{proposition}

Algorithms that efficiently find the block decomposition of a graph
are already known, see for instance \cite{HopTarjan} and
\cite[Chapter 4]{West}.
Let us briefly describe one possible approach to find
such decomposition, which we shall call the \emph{block
algorithm}. The description of the algorithm is first given in
the deterministic setting, and is later extended to the random
setting.

Motivated by \propref{prop:block_decomposition} (v),
we say that an unordered pair of vertices $uv$, where $u,v\in V(H)$, is
\emph{crossing} for $H$ if $u$ and $v$ lie in the same connected
component of $H$, and there is no block $B$ in $\block{H}$ containing
both $u$ and $v$. Note that by \propref{prop:block_decomposition} (ii),
a crossing pair is necessarily a non-edge of $H$. Another important
property of crossing pairs is the following.

\begin{proposition}

  \label{prop:crossing}
  Let $e$ be a crossing pair for $H$. Then the number of blocks of
  $H + \{e\}$ is strictly smaller than the number of blocks of $H$.

\end{proposition}

The input of the block algorithm consists of a pair $(H,H^*)$ of
graphs, where $H$ is a spanning subgraph of $H^*$ having the same number
of connected components as $H^*$. This requirement might seem rather
artificial at first, but it greatly simplifies the description of
the algorithm.
The output of the block algorithm is a graph $M$ such that
$H \subseteq M \subseteq H^*$ and $\block{M} \simeq \block{H^*}$. Moreover
$M$ is a minimal subgraph satisfying these properties, i.e., no proper subgraph of $M$
containing $H$ has the same number of blocks as $H^*$.

Let $M$ be the running graph. Initially we have $M := H$.
The main loop of the algorithm proceeds as follows.
\begin{quote}
If there exists a crossing pair $e\in E(H^*)\setminus E(M)$ for the
graph $M$, we add $e$ to $M$ and iterate the loop again.
Otherwise we stop and output $M$.
\end{quote}

Clearly, at the end of the algorithm we obtain a graph $M$ satisfying
the required properties. Moreover, by \propref{prop:crossing},
the number of iterations performed by the algorithm is less than the
number of blocks of $H$, as every new edge added to the running graph reduces the number
of blocks of the graph $M$.

In the random setting, the input parameter $H^*$ is a partially exposed
random graph, and $H$ is the graph containing the successfully
exposed edges from $H^*$. As we did in the DFS algorithm, we
only expose the edges of $H^*$ when their existence is queried by the
algorithm. One subtlety that should be remarked is that $\block{H^*}$ is
not known a priori, since the graph $H^*$ is random. The algorithm works
regardless. Moreover, \propref{prop:crossing} implies the following.

\begin{proposition}

  \label{prop:block_algorithm}
  The number of edges successfully tested by the block algorithm with
  input $(H, H^*)$ is less than the number of blocks of $H$.

\end{proposition}

Recall that we need to ensure that $H$ and $H^*$ have the same
number of connected components. To guarantee this assumption,
before we start the block algorithm, we run the DFS on $H^*$ and we always
choose an input parameter $H$ that contains the rooted spanning forest
produced by the DFS.

\subsection{P\'{o}sa's rotation-extension technique}
\label{sec:posa}

In this section we present yet another technique for showing the
existence of long paths and cycles in graphs. This technique was
introduced by P\'{o}sa~\cite{Posa} in his research on Hamiltonicity
of random graphs.

In quite informal terms, P\'{o}sa's lemma guarantees that expanding
graphs not only have long paths, but also provide a very convenient
structure for augmenting a graph to a Hamiltonian one by adding new
(random) edges.
To formalize this assertion, we need some definitions.
A graph $H$ is an \emph{$(m,2)$-expander} if $|N_H(X)| \ge 2|X|$
holds for every subset $X\subseteq V(G)$ of size $|X|\le m$.
Given a non-Hamiltonian graph $H$, a non-edge $e$ of $H$ is called
a \emph{booster} if $H+\{e\}$ is either Hamiltonian, or contains a path
which is longer than any path in $H$. The following consequence of
P\'{o}sa's technique (see, e.g., \cite[Lemma 8.5]{BolBook}) shows that every
connected and non-Hamiltonian graph $H$ with good expansion properties
has many boosters.

\begin{lemma}

  \label{lem:boosters}
  If $H$ is a connected non-Hamiltonian $(m,2)$-expander,
  then the number of boosters for $H$ is at least $(m+1)^2/2$.

\end{lemma}

\section{Proof of the main result}
\label{sec:proof}

For the rest of the paper, let $\eps=\eps(k):=\log^{-\frac{1}{10}} k$
and let $n$ be the number of vertices of $G$. We begin with the analysis
of the structure of $G$.
For that purpose, we make use of the following definition.

\begin{definition}

  \label{def:pseudoclique}
  A subset $C \subseteq V(G)$ of the vertices of $G$ is a
  \emph{pseudo-clique} if its size is bounded by $(1-4\eps)k<
  |C|\le (1+\eps)k$, and the minimum degree of $G[C]$ is at least
  $(1-4\eps)k$.

\end{definition}

This important notion plays a fundamental role in our analysis of $G$.
We later prove that if $G$ is not covered by many pseudo-cliques
with very few remaining vertices uncovered, then \aas $G_p$ contains
a cycle of length at least $k+1$.
To make this statement more precise, let $\mc C$ be a
collection of vertex-disjoint pseudo-cliques in $G$ such that
the union of their vertices $\bigcup_{C \in \mc C} C$
has maximum size. Vertices of $G$ not in $\bigcup_{C \in \mc C} C$
are called \emph{outcast} vertices, and let $\ell$ denote the number
of such vertices. We prove the following.

\begin{lemma}

  \label{lem:outcast_big}
  If $\ell > 10^7\cdot\frac{n}{\eps k}$ then \aas $G_p$ has a cycle of
  length at least $k+1$.

\end{lemma}

For the case when $\ell$ is small, we have the following.

\begin{lemma}

  \label{lem:outcast_small}
  If $\ell \le 10^7\cdot\frac{n}{\eps k}$, then either \aas $G_p$
  contains a cycle of length at least $k+1$, or there exist a
  pseudo-clique $C\in \mc C$ and a set $N$ of size $|N| \le 10$ such that
  there are at most $\eps k$ edges in $G$ connecting $C\setminus N$ to
  vertices not in $C\cup N$.

\end{lemma}

But if there exists such a pair $(C,N)$ as stated in
\lemref{lem:outcast_small}, $G_p$ must also have a cycle
of length at least $k+1$ a.a.s., as the next lemma shows.

\begin{lemma}

  \label{lem:good_pair}
  If there exist a pseudo-clique $C\in \mc C$ and a set $N\subseteq V(G)$
  of size at most $10$ such that $e_G\lt(C \setminus N, V(G)\setminus (C\cup N)\rt)
  \le \eps k$, then \aas $G_p[C\cup N]$ has a cycle of length at least $k+1$.

\end{lemma}

One can verify that lemmas \ref{lem:outcast_big}, \ref{lem:outcast_small},
and \ref{lem:good_pair} together imply \thmref{thm:main}. In the
next subsections, we devote ourselves to the proofs of these lemmas.
Our argument is divided into six steps. In each step, we may reveal
a portion of $G_p$ by testing some of the edges from $G$. The six
steps are:

\Step[Step 1] Pseudo-cliques were named for one clear reason: with
respect to $G_p$ they behave similarly as if they were cliques.
We formalize this claim by exposing the edges inside pseudo-cliques
and showing that a typical pseudo-clique contains a relatively long
cycle in $G_p$. We further delete from $G$ few vertices such that
in the remainder, every pseudo-clique induces a (large) Hamiltonian graph
in $G_p$. Finally, we prove that this deletion does not affect the
host graph much.

\Step[Step 2] We run a modified DFS algorithm on the resulting
graph from Step 1, handling pseudo-cliques as if they were
single vertices.
This way, the number of edges revealed in this step is small and
bounded by a function that depends only on $\ell$ and the number
of pseudo-cliques.

\Step[Step 3] We proceed with the block algorithm. The number of edges
that are revealed in this step is bounded similarly as in Step 2.
Hence after this step, we know the vertex sets of the blocks of $G_p$,
and \aas most outcast vertices still have almost $k$ untested
edges incident to them.

\Step[Step 4] The study of the internal structure of the blocks
provides some insight on how pseudo-cliques can interact with each other
and with other cycles. For instance, we prove that if a block contains
at least two pseudo-cliques then we already have exposed all the
edges of a cycle of length at least $k+1$.

\Step[Step 5] We use the results from the previous step, combined
with some double-counting arguments to prove Lemmas \ref{lem:outcast_big}
and \ref{lem:outcast_small}.
The only remaining case for the next step is the existence of a block
in our graph with one pseudo-clique, just a constant number of outcast
vertices, and only few edges between the pseudo-clique and the vertex set
outside the block.

\Step[Step 6] Finally, we analyze the case that remained after the
previous step. In some sense, this case is very close to the usual $\gnp(n,p)$
model: almost all vertices have degree close to $k$ inside the block,
and almost no edges leave the pseudo-clique to the outside of the
block. Using expansion properties of the random subgraph of the block,
we show that also in this case, we find a cycle of length at least
$k+1$ asymptotically almost surely.

\subsection{Step 1: preparing the pseudo-cliques}
\label{sec:step1}

Pseudo-cliques behave similarly as if they were cliques in $G$.
When exposed in $G_p$, pseudo-cliques typically contain large cycles of
length close to $k$. However, there might be a certain small
proportion of them behaving not in this typical way.
The aim of this subsection is to show that this seldom happens,
and therefore does not affect the remainder of the graph much.

Formally, let us consider a two-round exposure process.
Recall that we fixed a collection $\mc C$ of disjoint pseudo-cliques.
In the first round, we test edges inside pseudo-cliques with probability
$p_1$, where $p_1$ is such that $1-p=(1-p_1)^2$. Observe that
$p_1$ is roughly $\frac{\log k}{2k}$ and testing an edge with probability
$p$ (unsuccessfully) is the same as testing it twice (unsuccessfully)
with probability $p_1$.
Denote by $G^-$ the resulting random subgraph.
Let $W_1$ be the set of vertices that have degree at most
$\log k/100$ inside their pseudo-cliques in $G^-$.

In the second round we again expose with probability $p_1$ the edges
inside pseudo-cliques in $\mc C$ that were not successfully exposed during
the first round;
the resulting supergraph of $G^-$ is denoted by $G^+$. For technical
reasons, we would like the remainders $C\setminus W_1$ of pseudo-cliques
$C\in \mc C$ to satisfy the properties:
\begin{enumerate}[\bf \hspace{5pt} (P1)] \smalllist
\item $C\cap W_1$ has fewer than $\eps k/2$ vertices,
\item $e_{G^-}(X,Y) > 0$ for any two disjoint sets
$X,Y\subseteq C\setminus W_1$ of size at least $6\eps k$,
\item the induced graph $G^+\lt[C\setminus W_1\rt]$ is Hamiltonian.
\end{enumerate}
We now define the set $W_2$ to be the union of those pseudo-cliques
$C\in \mc C$, for which the above properties do not simultaneously
hold for $C\setminus W_1$. We refer to the set $W:=W_1\cup W_2$ as the
\emph{waste}. The set $W$ contains the vertices we aim to delete from the graph
$G$ to obtain the new graph $G':=G\setminus W$. Finally, let $Z_1$ be the
set of all outcast vertices $u$ such that at least an
$\frac{\eps}{3}$-proportion of its neighbors from $G$ belong to $W$.
The probability that $u \in Z_1$ is bounded by the following statement.

\begin{lemma}
\label{lem:prob_z1}
Let $u$ be an outcast vertex. Then $\Prob[u\in Z_1] \le 1/k^3$.
\end{lemma}

We split the proof of \lemref{lem:prob_z1} into several propositions,
from which the statement of the lemma is a trivial consequence.
The first proposition of the series insures that \aas most outcast vertices
do not have many neighbors in $W_1$.

\begin{proposition}

  \label{prop:high_degree_waste}
  Let $u$ be an outcast vertex and denote by $d\geq k$ its degree in
  $G$. The probability that at least $\eps d/6$ neighbors of $u$ belong
  to $W_1$ is at most $1/k^4$.

\end{proposition}
\begin{proof}
The probability that a vertex $v$ from a pseudo-clique $C$ has degree
at most $\log k/100$ in $G^-[C]$ is already sufficiently small.
However, these events are not independent: the event that $v$ has
small degree in $G^-[C]$ is positively correlated with another vertex
from the same pseudo-clique getting small degree in $G^-[C]$.
Since the statement of the proposition is far from being tight,
one possibility to overcome this technicality is the following.
Let $\vec{G}$ be the digraph obtained from $G$ by replacing each edge
$vw\in E(G)$ with two oriented edges $\vec{vw},\vec{wv}\in E(\vec{G})$.
We test each of the $2\lt|E\lt(G[C]\rt)\rt|$ oriented edges corresponding to the
edges of $G[C]$ independently with probability $p_2$, where $p_2$ is such
that $1-p_1 = (1-p_2)^2$, and roughly $p_2 \approx \frac{p_1}{2}
\approx \frac{\log k}{4k}$.
Next, we say that we successfully exposed the (non-oriented) edge $vw\in E(G)$
if we successfully exposed at least one of the oriented edges
$\vec{vw}$ or $\vec{wv}$. In this model, all non-oriented edges are exposed
independently at random with probability $p_1$.
Thus, we can assume that each edge $vw$ of $G[C]$ that became a non-edge
also had two corresponding oriented non-edges, $\vec{vw}$ and $\vec{wv}$,
in the random digraph.
Hence, in order for $v$ to get at most $\log k/100$ non-oriented edges,
all but at most $\log k/100$ of the oriented edges going out from $v$
to other vertices of $C$ must become non-edges.
Now, these events (``all but at most $\log k/100$ oriented edges going
out from a fixed vertex from $C$ to other vertices in $C$ were tested
as non-edges'') are indeed independent for any two vertices from $C$.

For one vertex $v\in C$, since the minimum degree in $G[C]$ is at least
$(1-4\eps)k$, the probability of this event is at most
\begin{equation}
  \label{eqn:prob_vtx_waste1}
  \Prob\lt[\Bin((1-4\eps)k, p_2) \le \log k / 100\rt] < k^{-1/5}
\end{equation}
due to \lemref{lem:chernoff}. Thus, the probability that at least
$\eps d / 6$ neighbors of $u$ belong to $W_1$ is bounded by
$\Prob\lt[\Bin(d, k^{-1/5}) > \eps d / 6 \rt]$, and another
application of \lemref{lem:chernoff} finishes the proof of
the proposition.
\end{proof}

For a pseudo-clique $C$, let us denote by $C^-$ the remainder
$C\setminus W_1$. Similarly to \propref{prop:high_degree_waste},
we need to ensure that also for a vertex from a pseudo-clique $C$,
after the first round of exposure, \aas only few neighbors of this
vertex are in $C\cap W_1$. The proof of this proposition follows
the lines of the proof of \propref{prop:high_degree_waste} and is
therefore omitted.

\begin{proposition}

  \label{prop:min_degree_nonwaste}
  For fixed $C\in \mc C$ and $u\in C$, the probability
  that in $G^-$, at least $\log k/200$ of the neighbors of $u$ are in
  $C\cap W_1$, is at most $1/k^6$.

\end{proposition}

We remark that one could have replaced $\log k /200$ by a large constant in
the statement of \propref{prop:min_degree_nonwaste}. Indeed, the number of
neighbor of $u$ in $C\cap W_1$ can be roughly bounded by a binomial
random variable of $O(\log k)$ trials with success probability
$k^{-1/5}$. However, we do not require such tight estimates.

The very same calculation also shows that \aas $C^-$ is large enough,
as required to satisfy (P1).

\begin{proposition}

  \label{prop:size_C_minus}
  For fixed $C\in \mc C$, with probability at least $1-1/k^5$,
  we have $|C \cap W_1|< \eps k/2$.

\end{proposition}

Notice that the inequalities in \propref{prop:size_C_minus} are
again far from being sharp, but they already suffice for our purposes.

Recall that for the second property (P2), we need
$G^-[C^-]$ to have edges between any two reasonably large
disjoint sets. The next proposition ensures that \aas this is indeed
the case.

\begin{proposition}

  \label{prop:property_p2}
  For every $C\in \mc C$, with probability at least $1-1/k^5$,
  we have $e_{G^-}(X,Y) > 0$ for any two disjoint sets $X,Y \subseteq C^-$,
  each of size at least $6 \eps k$.

\end{proposition}
\begin{proof}
In $G$, for every choice of the sets $X,Y \subseteq C$, we have
$e_G(X,Y) \ge 6 \eps^2 k^2$, as every vertex from $X$ has at least
$|Y| - 5 \eps k \ge \eps k$ neighbors in $Y$.
This is because every vertex in a pseudo-clique $C$ has at most
$5\eps k$ non-neighbors in $G[C]$.
Thus, the probability that $e_{G^-}(X,Y) = 0$ is at most
$(1-p_1)^{6\eps^2k^2}\le \exp(-k \sqrt{\log k})$.
Since there are at most $4^k$ possible choices for the pair
$X,Y$, a simple application of the union bound finishes the proof.
\end{proof}

For the last property (P3), required to ensure that a pseudo-clique $C$ is
not put into $W_2$, we need $G^+[C^-]$ to be Hamiltonian. To prove
the Hamiltonicity of $G^+[C^-]$ we first show in the next proposition
that $G^-[C^-]$ is a good expander.

\begin{proposition}

  \label{prop:expander}
  For $C \in \mc C$, with probability at least $1-3/k^5$,
  the induced graph $G^-[C^-]$ is a $(k/6000,2)$-expander.

\end{proposition}
\begin{proof}
Suppose that there exists a set $A\subset C^{-}$ of size
$|A|\le k/6000$ such that $\lt|N_{G^-[C^-]}(A)\rt|< 2|A|$.
Also assume that the conclusion of \propref{prop:min_degree_nonwaste}
does not hold for any vertex in $C^-$, i.e., no vertex in $C^-$ has
more than $\log k /200$ neighbors in $C\cap W_1$.
This happens with probability at least $1 - 2/k^5$ by the union bound.
Thus, if $u\in C^{-}$, we have $\deg_{G^-[C^-]}(u) \ge \deg_{G^{-}}(u)
-\log k /200 \ge \log k / 200$.
Now let $B = A \cup N_{G^-[C^-]}(A)$. Then $|B| < 3|A|
\le k / 2000$, and $\lt|E\lt(G^-[B]\rt)\rt|\ge|A| \log k /400\ge
|B| \log k/1200$. On the other hand, by \lemref{lem:chernoff2} and
the union bound, we have
\begin{align*}
  &\Prob\lt[\exists B\subset C,~|B|\le k/2000:
    ~\lt|E\lt(G^-[B]\rt)\rt|\geq |B| \log k/1200\rt]\\
  &\qquad
    \le \sum_{b \le k/2000} \binom{|C|}{b}\Prob\lt[
    \Bin\lt(\binom{b}{2}, p_1\rt) > b \log k /1200\rt] \\
  &\qquad
   \le \sum_{b\le k/2000} \binom{|C|}{b}\lt(\frac{e p_1\binom{b}{2}}{b
   \log k/1200}\rt)^{b\log k/1200} < 1/k^5,
\end{align*}
where in the last inequality we used that $|C| \le (1+\eps)k$
and that $p_1\approx \frac{\log k}{2k}$. This concludes the proof
of the proposition.
\end{proof}

Finally, we show that with sufficiently high probability, $G^+[C^-]$
is Hamiltonian. Notice that we could strengthen the statement
and ask for $G^+[C^-]$ to be Hamilton connected. However, Hamiltonicity
suffices for our proof, and it is technically slightly easier to show.

\begin{proposition}

  \label{prop:good_clique}
  For every $C\in \mc C$, with probability at least $1-5/k^5$
  all properties {\rm (P1)}, {\rm (P2)}, and {\rm (P3)}
  hold for $C$.

\end{proposition}
\begin{proof}
After propositions \ref{prop:size_C_minus}, \ref{prop:property_p2},
and \ref{prop:expander}, we can assume that $G^-[C^-]$ is a connected
$(k/6000,2)$-expander on at least $(1-5\eps)k$ vertices, and satisfies
properties (P1) and (P2). The connectivity of $G^-[C^-]$ is a consequence
of \propref{prop:expander}, which implies that every connected component
of $G^-[C^-]$ has at least $k/2000$ vertices, together with
\propref{prop:property_p2}.  Conditioned on these assumptions,
we would like to show that then $G^+[C^-]$ is Hamiltonian with
probability at least $1-1/k^5$. Indeed, in case a supergraph $H$ of
$G^-[C^-]$ is not Hamiltonian, \lemref{lem:boosters} guarantees a
quadratic number of boosters.
Now, let us look at the second round of exposure as a random process,
with non-edges of $G^-[C^-]$ turning into edges one-by-one, analogous
to the standard random process coupling $\gnp(n,p)$ and $\gnp(n,M)$.
The new edges are exposed in a random order, their number
$\lt|E\lt(G^+\lt[C^-\rt]\rt)\setminus E\lt(G^-\lt[C^-\rt]\rt)\rt|$
is binomially distributed, thus by \lemref{lem:chernoff} with
probability at least $1-e^{-k}$, there are $\Omega(k \log k)$
new successfully exposed edges.
After every exposed edge, we update the set of boosters ---
keeping in mind that there are still quadratically many of them.
Hence, every successfully exposed edge is a booster with probability
at least a constant bounded away from zero.
Thus we expect that the number of additional exposed edges needed for the
graph induced by $C^-$ to become Hamiltonian is at most linear.
Furthermore, we can use \lemref{lem:chernoff} to say that the
probability that we expose $\omega(k)$ edges and we do not make the
graph on $C^-$ Hamiltonian, is at most $e^{-k}$, and the statement
of the proposition follows.
\end{proof}

The following statement can be derived in the same way
as \propref{prop:good_clique}, hence we omit its proof.

\begin{proposition}

  \label{prop:small_clique}
  We may assume that there is no set
  $X\subseteq V(G)$ of size $(1+\eps/2)k \le |X| \le
  (1+20000\eps)k$ such that the minimum degree of $G[X]$ is at least
  $(1-10\eps)k$, as otherwise $G_p[X]$ \aas would contain a cycle
  of length at least $k+1$.

\end{proposition}

The last proposition allows us to further assume from this point on
that all pseudo-cliques in $\mc C$ have size less than $(1+\eps/2)k$.
We are ready to prove \lemref{lem:prob_z1}.

\begin{proof}[Proof of \lemref{lem:prob_z1}{}]
Let $u \in Z_1$, and let $d$ denote the degree of $u$ in $G$.
Our aim is to bound the number of neighbors of $u$ that are in $W$.
We remark that the following estimations for the number of neighbors
of $u$ which belong to $W$ are true even if we drop the assumption
that $u$ is outcast.

Either $u$ has $\eps d / 6$ neighbors in $W_1$, or it has the same amount of
neighbors in $W_2$. \propref{prop:high_degree_waste} bounds the probability
of the first case to happen by at most $1/k^4$. For the second case,
notice that \propref{prop:good_clique} implies that $\Prob\lt[w \in W_2\rt]
\le \frac{5}{k^5}$ for all $w\in \bigcup_{C\in \mc C} C$. By Markov's inequality,
we have that $\Prob\lt[\lt|N(u) \cap W_2\rt| > \eps d / 6\rt] < \frac{30}{\eps k^5}
< \frac{1}{k^4}$. Therefore, by the union bound, $\Prob[u\in Z_1] < 1/k^3$,
concluding the proof of the lemma.
\end{proof}

\lemref{lem:prob_z1} bounds the number of outcast vertices that lost
a significant proportion of their neighbors after the deletion
of the waste from $G$ to obtain $G'$. By Markov's inequality,
asymptotically almost surely, the size of $Z_1$ is bounded by
\begin{equation}
  \label{eqn:size_z1}
  |Z_1| \le \frac{n}{k^2}.
\end{equation}
This inequality tells us that the influence of the waste is not
too large, so for most of our subsequent arguments, we
can completely ignore the vertices from $W$. Also, from the
definition of $Z_1$, if $v\in G'\setminus Z_1$ is an outcast
vertex then $\deg_{G'}(v)\ge \left(1-\frac{\eps}{3}\right)k$,
hence $v$ still retains most of its degree after the deletion of $W$.
However, in the final part of the proof of our main theorem,
we have to use the full structure of $G$ and incorporate the waste
vertices back. Therefore, we need a lemma to state what typically
happens to a pseudo-clique after we delete the vertices from the
waste.

\begin{lemma}

  \label{lem:edgesoutside}
  Consider an arbitrary pseudo-clique $C\in \mc C$, and denote by $D_1$ the set
  of vertices from $C$ having more than $\eps k$ neighbors in $G$ outside
  of $C$. Let $D'_1$  be the union of $D_1\cap W$ together with the vertices
  in $D_1\setminus W$ that lost more than a $\frac{1}{100}$-proportion of its
  neighbors outside of $C$ after the removal of the waste vertices from $W$.
  Furthermore, let $D_2$ be the set of vertices not in $C$ that have
  at least $\eps k$ neighbors in $C$ in the graph $G$. Finally,
  let $\mc E$ denote the set of edges from $G$ connecting
  $C\setminus D_1$ to a vertex not in $C \cup D_2$. Then \aas we have
  \begin{equation}
    \label{eqn:edgesoutside}
    |\mc E \setminus E(G')| \le |\mc E|/100, \quad |D_2\cap W|\le
    |D_2|/100 \quad\text{and}\quad |D'_1| \le |D_1|/100.
  \end{equation}
  Therefore, \aas at least $(1+o(1))|\mc C|$ pseudo-cliques in $\mc C$ satisfy
  \eqref{eqn:edgesoutside}.

\end{lemma}
\begin{proof}[Sketch of the proof]
For each fixed vertex $u$, the probability that $W$ contains $u$ is
either zero (if $u$ is outcast) or tiny, as the inequality
\eqref{eqn:prob_vtx_waste1} together with \propref{prop:good_clique}
imply that both $\Prob[u\in W_1]$ and $\Prob[u\in W_2]$ are small.
Similarly, for each fixed edge $e$, the probability that one of
its endpoints belongs to $W$ is also very small. In expectation,
we have $\E[\mc E\setminus E(G')] = o(|\mc E|)$, hence by Markov's
inequality we know that \aas $|\mc E\setminus E(G')| \le |\mc E|/100$.
Similarly, we have $\E[D_1\cap W]= o(|D_1|)$ and $\E[D_2\cap W]=o(|D_2|)$.
Moreover, for each vertex $u\in D_1$, if we denote by $d$ the number of
edges connecting $u$ to a vertex outside of $C$, and by $d'$ the number
of edges connecting $u$ to a vertex in $W\setminus C$, then
$\E[d'] = o(d)$. Thus, by Markov's inequality, we know that
$\Prob[v\in D'_1] = o(1)$, and the lemma follows by another application
of Markov's inequality.
\end{proof}

\subsection{Step 2: exploring the connected components}
\label{sec:step2}

Recall that at this point, some of the edges of $G$ were already
tested in $G_p$, namely all the edges inside pseudo-cliques from
$\mc C$. Let $Q_1$ be the set of tested edges from the partially exposed
$G_p$ that live inside $G'$, and let $E_1 \subseteq Q_1$ be the subset
of the successfully tested edges. To find the connected components of
the partially exposed $G'_p=G_p\setminus W$ using the DFS algorithm,
we adopt the following DFS priority:
\begin{quote}
Whenever the DFS reaches a vertex $v$ from a pseudo-clique
$C\in \mc C$, the algorithm, instead of testing new edges, walks
through an already exposed Hamilton cycle in $G_p[C\setminus W]$
using the edges from $E_1$, until it visits all vertices from
$C \setminus W$.
\end{quote}
In the rooted spanning output forest $T$, this Hamilton cycle forms a path,
and the algorithm saved many edge tests this way. This observation is
stated more formally as follows.

\begin{observation}

  \label{obs:pseudo_path}
  For every pseudo-clique $C\in \mc C$ such that $C\not\subseteq W$,
  there exists a path in $T$ whose vertices are precisely the
  vertices in $C\setminus W$.

\end{observation}

Let $Q_2'$ be the set of tested edges, and let $E_2 \subseteq Q_2'$
be the set of successfully tested edges in this exploration of $G'_p$
by the DFS. Clearly $E(T)\subseteq E_1\cup E_2$, and $|E_2| < \ell + |\mc C|$,
since once we reach a pseudo-clique $C$, we do not need to test edges
until all the vertices of $C \setminus W$ are reached.

Next, we query all the untested edges connecting vertices from $G'$ which
have distance at least $k + 1$ with respect to the forest $T$. Let $Q_2''$
be the set of all such edges. We test the edges in $Q_2''$ one by one,
in an arbitrary order. If by chance we successfully expose one edge from
$Q_2''$, we automatically obtain a cycle of length at least $k+1$ in
$G'_p$, as desired in \thmref{thm:main}, and we stop the whole procedure.
In particular, the total number of edges in $Q_2''$ must be very small,
say $|Q_2''| < \eps k$, as otherwise we \aas would have a long cycle.
Let $Q_2$ be the union of $Q_2'$ with the tested edges from
$Q_2''$. We can estimate the total number of edges $Q_2$ using
\propref{prop:dfs_new_edge} and \lemref{lem:chernoff}, obtaining the next
statement.

\begin{corollary}

  \label{cor:tested_dfs}
  Asymptotically almost surely, we have
  \[
    |Q_2| \le \frac{1.1}{p}\cdot |E_2| + \eps k
    \le  \frac{1.2}{p}\lt(\ell + \frac{n}{k}\rt) + \eps k.
  \]
  Moreover, if not all the edges in $Q_2''$ were tested
  at this point, then we already have exposed a cycle
  of length at least $k+1$ in $G'_p$.

\end{corollary}

We would like to remark that the expression $\ell + \frac{n}{k}$ is not
guaranteed to tend  to infinity with $k$, so the inequality $|Q_2'| \le
\frac{1.1}{p}\cdot |E_2|$ is not guaranteed to hold asymptotically almost
surely. However, \corref{cor:tested_dfs} is true because of the extra
$\eps k$ term, as
\[
 \Prob\lt[\Bin\lt(\frac{1.2}{p}\cdot \lt(\ell +
 \frac{n}{k}\rt) + \eps k, p\rt) \ge \ell + \frac{n}{k}\rt] \to 1
 \quad\text{as}\quad k\to\infty.
\]

\subsection{Step 3: the block decomposition}
\label{sec:step3}

In this subsection, we apply the block algorithm to the input $(H,H^*)$,
where $H=T+E_1$ and $H^*$ is the current partially exposed $G'_p$. Recall
that $T+E_1$ might have some large cycles already, coming from the
exposed pseudo-cliques in Step 1. Thus we can bound the number of blocks
of $T+E_1$ from above by $\ell + |\mc C|$. This is because for every
pseudo-clique $C \in \mc C$ which is not completely inside the waste
$W$, the vertex set $C\setminus W$ necessarily induces a Hamiltonian graph
$(T+E_1)[C\setminus W]$.

Let $Q_3$ be the set of edges from $G'_p$ tested during the execution
of the block algorithm, and let $E_3\subseteq Q_3$ be the subset of
the successfully tested edges. From \propref{prop:block_algorithm},
we know that the number of successfully tested edges revealed by
the block algorithm is at most the number of blocks of $T+E_1$. Moreover,
by the observation discussed in the last paragraph, we also know
that the total number of blocks of $T+E_1$ is at most $\ell + |\mc C|$,
hence $|E_3| \le \ell + |\mc C|$, and by \lemref{lem:chernoff} we have
the following corollary.

\begin{corollary}

  \label{cor:tested_block}
  Asymptotically almost surely $|Q_3| \le \frac{1.2}{p}\left(\ell +
  \frac{n}{k}\right) + \eps k$.

\end{corollary}
We would like to draw the reader's attention to the fact that we added the
term $\eps k$ to the right hand side of the inequality in \corref{cor:tested_block}.
This is because we want to make sure that the right side, when multiplied
by $p$, tends to infinity with $k$. We recall that a similar ``trick'' was
used in \corref{cor:tested_dfs}.

\subsection{Step 4: the structure inside the blocks}
\label{sec:step4}

Let $\mc B$ be the family of all blocks of the partially exposed $G'_p$
obtained in Step 3. Here, the edges of every block $B\in \mc B$ consist
of those successfully exposed in $G'_p$ so far, i.e., $E(B)\subseteq E_1
\cup E_2 \cup E_3$. Moreover, the cut-vertices of $G'_p$ are precisely the
same as the cut-vertices of $T+(E_1\cup E_3)$. One should also
observe the following.

\begin{observation}

  \label{obs:block_properties}
   For each $B\in \mc B$, the graph $T\cap B$ is a tree.
   Moreover, if $B_1$ and $B_2$ are two distinct blocks from $\mc B$
   having a vertex $v$ in common, then $v$ is the smallest vertex
   (with respect to $\le_T$) from at least one of the two blocks
   $B_1$ or $B_2$.

\end{observation}

\begin{figure}[ht!]
\centering
\begin{tikzpicture}
  [scale=0.8, auto=left, every node/.style={circle, draw, fill=black,
   inner sep=0pt, minimum width=4pt}]
  \node (n1_1) at (-0.2,  0) {};
  \node (n1_2) at (0.2, 0) {};
  \node (n2) at (-1.2, -1) {};
  \node (n3) at (1.2, -1) {};
  \node (n4) at (-2.2, -2) {};
  \node (n5_1) at (-0.2, -2) {};
  \node (n5_2) at (-0.4, -2.2) {};
  \node (n5_3) at (0.0, -2.2) {};
  \node (n6) at (-1.4, -3.2) {};
  \node (n7) at (0.0, -3.2) {};
  \node (n8) at (-0.5, -4.2) {};
  \node (n9) at (0.5, -4.2) {};
  \node (n10_1) at (2.2, -2) {};
  \node (n10_2) at (2.2, -2.4) {};
  \node (n11) at (2.2, -3.4) {};

  \foreach \from/\to in {n1_1/n2, n1_2/n3, n2/n4, n2/n5_1,
                         n5_2/n6, n5_3/n7, n7/n8, n7/n9,
                         n3/n10_1, n10_2/n11}
    \draw (\from) -- (\to);

  \draw[style=dashed] (0, 0) circle (0.5);
  \draw[style=dashed] (-0.2, -2.1) circle (0.5);
  \draw[style=dashed] (2.2, -2.2) circle (0.5);

  \node [draw=none,fill=none,label=above:root] at (0,0.2) {};

  \node [label=above:root] (m1) at (6, 0) {};
  \node (m2_1) at (6, -1) {};
  \node (m2_2) at (5.8, -1.2) {};
  \node (m2_3) at (6, -1.4) {};
  \node (m2_4) at (6.2, -1.2) {};
  \node (m3) at (4.8, -2.2) {};
  \node (m4) at (6, -2.4) {};
  \node (m5) at (7.2, -2.2) {};
  \node (m6) at (5.5, -3.4) {};
  \node (m7) at (6.5, -3.4) {};

  \foreach \from/\to in {m1/m2_1, m2_2/m3, m2_3/m4, m2_4/m5,
                         m4/m6, m4/m7}
    \draw (\from) -- (\to);

  \draw[style=dashed] (6, -1.2) circle (0.5);

\end{tikzpicture}
\caption{Blocks and cut-vertices of $G'_p$ together with the
rooted forest $T$.}
\label{fig:blocks}
\end{figure}
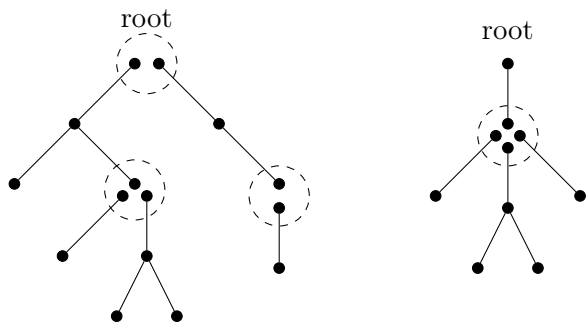

The content of the previous observation is illustrated in
\figref{fig:blocks}. In the picture, each connected component
represents a subtree of the form $T\cap B$ for some $B\in \mc B$.
The dashed ovals represent the cut-vertices from $G'_p$ (all small
solid circles inside the dashed ovals actually represent the same
cut-vertex).

Note that for all pseudo-cliques $C \in \mc C$ with
$C \not\subseteq W$, there exists a unique block $B \in \mc B$ such
that $C\setminus W\subseteq V(B)$. This is because cycles are $2$-connected.
In this case, with slight abuse of notation, we say that $B$ \emph{contains}
the pseudo-clique $C$.

The next proposition shows that a block containing
more than one pseudo-clique already has a long cycle.

\begin{proposition}

  \label{prop:one_block_one_clique}
  If $B \in \mc B$ is a block that contains two distinct
  pseudo-cliques $C_1,C_2\in \mc C$, then $B$ contains a cycle of
  length at least $k+1$.

\end{proposition}

Before we prove \propref{prop:one_block_one_clique}, let us
prove an auxiliary statement.

\begin{proposition}

  \label{prop:rotating_cliques}
  Let $C \in \mc C$ be such that $C\not\subseteq W$, and let
  $B\in \mc B$ be the block containing $C$.
  Then for every two distinct vertices $u,v\in C \setminus W$,
  the induced graph $B[C \setminus W]$ contains a path of length at least
  $(1-20\eps)k$ connecting $u$ to $v$.

\end{proposition}
\begin{proof}
We want to show that there exists a path $P$ in $B[C\setminus W]$
connecting $u$ and $v$ of length at least $(1-20\eps)k$. Because of
property (P3) stated in \secref{sec:step1}, we know that $B[C \setminus W_1]=
B[C\setminus W]$ is Hamiltonian. Let $J$ be a Hamilton cycle in
$C\setminus W$. Next, consider the two paths $P_1$ and $P_2$
obtained from the cycle $J$ connecting the vertices $u$ and $v$. Assume
that $P_2$ is no longer than $P_1$. By property (P1), $J$ is of length
at least $(1-5\eps)k$, so $P_1$ has at least $(1-5\eps)k/2$ vertices.
If the length of $P_1$ is greater than $(1-20\eps)k$, our
proposition immediately follows by taking $P:=P_1$.
Otherwise the path $P_2$ has at least $15\eps k$ vertices.
Let $X$ be the set of the $6\eps k$ vertices from the path $P_1$ which
are closest to the endpoint $v$. Similarly, let $Y$ be
the set of the $6\eps k$ vertices from the path $P_2$
which are closest to the other endpoint $u$. Using property (P2),
we know that $B[C \setminus W]$
has an edge $e=xy$ connecting a vertex $x$ from $X$ to a vertex
$y$ from $Y$, as shown in \figref{fig:rotation_pseudocliques}.

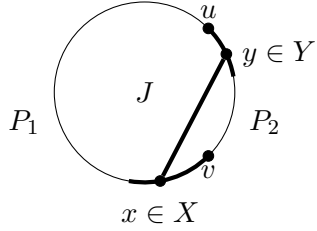
\begin{figure}[ht!]
\centering
\begin{tikzpicture}
  [scale=0.8, auto=left, every node/.style={circle, draw, fill=black,
   inner sep=0pt, minimum width=4pt}]
  \node[label=above:$u$] at ({1.5*cos(45)}, {1.5*sin(45)}) {};
  \node[label=below:$v$] at ({1.5*cos(-45)}, {1.5*sin(-45)}) {};
  \node[label={[label distance = -0.2cm]-90:$x\in X$}]
     (node_x) at ({1.5*cos(-80)}, {1.5*sin(-80)}) {};
  \node[label={[label distance = 0.1cm]0:$y\in Y$}]
     (node_y) at ({1.5*cos(25)}, {1.5*sin(25)}) {};

  \node[draw=none,fill=none] at (0, 0) {$J$};
  \node[draw=none,fill=none] at (-2, -0.5) {$P_1$};
  \node[draw=none,fill=none] at (2, -0.5) {$P_2$};

  \draw (0cm, 0cm) circle(1.5cm);
  \draw[ultra thick] (node_x) -- (node_y);

  \draw[ultra thick] ([shift=(45:1.5cm)]0,0) arc (45:10:1.5cm);
  \draw[ultra thick] ([shift=(-45:1.5cm)]0,0) arc (-45:-100:1.5cm);

\end{tikzpicture}
\caption{Rotation inside pseudo-cliques.}
\label{fig:rotation_pseudocliques}
\end{figure}

We can build the longer path $P$ by patching
two segments from $P_1$ and $P_2$ together with the edge
$e$ as follows. The initial segment of $P$ consists
of the path in $P_1$ connecting $u$ to $x$, while
the final segment of $P$ consists of the path in $P_2$
connecting $y$ to $v$, and these two segments are interconnected
by $e$. The total length of $P$ is at least the length of $J$ minus
$12 \eps k$, hence the length of $P$ is at least $|J|-12\eps k>
(1-20\eps)k$, finishing the proof of the proposition.
\end{proof}
We are ready to prove \propref{prop:one_block_one_clique}.
\begin{proof}[Proof of \propref{prop:one_block_one_clique}{}]
In the proof of this proposition, we use \thmref{thm:menger}
to merge long cycles. Since $B$ is $2$-connected,
\thmref{thm:menger} asserts the existence of two vertex
disjoint paths $P_1,P_2$ in $B$ connecting
$C_1 \setminus W$ to $C_2 \setminus W$. Let $u_1,v_1$ be
the endpoints of $P_1,P_2$ (respectively) in $C_1\setminus W$.
Similarly, let $u_2, v_2$ be the endpoints of $P_1,P_2$ in
$C_2\setminus W$.

By \propref{prop:rotating_cliques}, we can obtain two paths $P_3,P_4$
both having length at least $(1-20\eps)k$, where $P_3$ is a path
in $C_1\setminus W$ connecting $u_1$ to $v_1$, and $P_4$ is a path
in $C_2\setminus W$ connecting $u_2$ to $v_2$. By patching together
$P_1$, $P_3$, $P_2$, and $P_4$ in that order, we obtain a cycle of
total length larger than $(2-40\eps)k > k+1$, thereby proving the
proposition.
\end{proof}

As we have seen in the proof \propref{prop:one_block_one_clique},
we can use pseudo-cliques to obtain long cycles, which then can be
merged into even longer cycles. We do not need to use the full strength
of pseudo-cliques in order to merge cycles. In the proof of
\propref{prop:rotating_cliques}, the edge $e$ played an important role,
as it allowed us to ``rotate'' inside the relatively long cycle.
In what follows, we describe a weaker structure that also allows
this ``rotation'' operation. We say that a cycle $J$, formed by some of
the successfully exposed edges from a partially exposed $G'_p$, is a
\emph{rotating cycle} if all properties below hold simultaneously:
\begin{enumerate}[\bf \hspace{5pt} (P1$\star$)] \smalllist
\item $J$ has at least $(1-4\eps)k$, but at most $k$ vertices,
\item all but one edge of $J$ belong to the forest $T$ revealed
in Step 2,
\item if $u\in V(J)$ is the largest vertex with respect to the
order $\le_T$ (we call $u$ the \emph{pivot} of $J$), then there exists
at least $(1-4\eps) k$ untested edges in the partially exposed $G'_p$
connecting $u$ to another vertex of $J$.
\end{enumerate}
The properties listed previously bear some resemblance to the ones
enumerated in \secref{sec:step1}. For instance (P1) and (P1$\star$) both
state some bounds about the size of the structure under consideration.
Property (P2$\star$) might look somewhat artificial at first, but we
observe that for every pseudo-clique $C\in\mc C$ such that $C\not\subseteq W$,
the graph $G'_p[C \setminus W]$ contains a Hamilton path that
is entirely contained in $T$, a consequence of the priority of the
DFS remarked in \secref{sec:step2}.
Finally, (P3$\star$) is the property that will allow us to perform the
rotation per se, and note that (P3$\star$) clearly implies the lower bound
of the length of $J$ in (P1$\star$).

The next proposition describes the operation of rotation, which is
similar to the rotation described in \propref{prop:rotating_cliques}.

\begin{proposition}

  \label{prop:rotating_cycles}
  Let $J$ be a rotating cycle with pivot $u\in V(J)$, and fix any two
  distinct vertices $x,y \in V(J) \setminus \{u\}$. After exposing the
  untested edges connecting $u$ to the other vertices in $J$,
  \aas we can find a path in $G'_p[V(J)]$ between $x$ and $y$ of
  length at least $(2-10\eps)k/3$.

\end{proposition}
\begin{proof}
Let $P_1$ and $P_2$ be the two paths between $x$ and $y$ obtained
from the cycle $J$, with lengths $l_1$ and $l_2$, respectively.
Assume, without loss of generality, that $u \in V(P_1)$,
and that the distance from $u$ to $y$ is no larger than the distance
from $u$ to $x$ in the path $P_1$. Let $N$ be the set of all
vertices $w$ of $J$ such that $uw$ is an untested edge of $G'_p$,
and let $N_1=N\cap P_1$ and $N_2=N \cap P_2$.
If either $N_1$ or $N_2$ has size at least $(2-10\eps)k/3$
then we are done, since $l_1 \ge |N_1|$ and $l_2 \ge |N_2|$.
Otherwise, both $N_1$ and $N_2$ have size at least
$(1-2\eps)k/3$, because (P3$\star$) implies that
$|N| \ge (1-4\eps)k$. Next, we test all the edges connecting
$u$ to the $\eps k$ vertices from $N_2$ which
are closest to $x$ with respect to the path $P_2$. This is possible
because $|N_2| > \eps k$. Asymptotically almost surely, we can find
a successfully exposed edge $e=uw$ where $w$ belongs to this subset
of $N_2$ of size $\eps k$.

\begin{figure}[ht!]
\centering
\begin{tikzpicture}
  [scale=0.8, auto=left, every node/.style={circle, draw, fill=black,
   inner sep=0pt, minimum width=4pt}]
  \node[label=above:$x$] at ({1.5*cos(120)}, {1.5*sin(120)}) {};
  \node[label=below:$y$] at ({1.5*cos(-120)}, {1.5*sin(-120)}) {};
  \node[label=left:$u$] (node_u) at ({1.5*cos(190)}, {1.5*sin(190)}) {};
  \node[label=above:$w$] (node_w) at ({1.5*cos(100)}, {1.5*sin(100)}) {};

  \draw (0cm, 0cm) circle(1.5cm);
  \draw[ultra thick] (node_u) -- (node_w);

  \draw[ultra thick] ([shift=(120:1.5cm)]0,0) arc (120:190:1.5cm);
  \draw[ultra thick] ([shift=(-120:1.5cm)]0,0) arc (-120:100:1.5cm);

\end{tikzpicture}
\caption{Rotation inside the rotating cycle.}
\label{fig:rotating_cycle}
\end{figure}
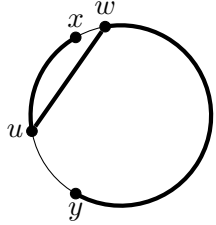

We can then obtain a path $P$ as follows: we use the segment
from $P_1$ connecting $x$ to $u$, and then we traverse the edge $e$,
to reach the vertex $w$, and then use the segment from $P_2$ connecting
$w$ to $y$ as illustrated in \figref{fig:rotating_cycle}.
The length of $P$ is at least $l_1/2 + |N_2| - \eps k \ge |N_1|/2+
|N_2|-\eps k=(|N|+|N_2|)/2-\eps k\ge(2 - 10\eps)k/3$, concluding
the proof.
\end{proof}

Analogous to pseudo-cliques, rotating cycles are also somewhat tied
to the block structure of $G'_p$. For every rotating cycle $J$,
there exists a unique block $B\in \mc B$ such that $V(J) \subseteq V(B)$.
The equivalent of \propref{prop:one_block_one_clique} for rotating
cycles is the next statement.

\begin{proposition}

  \label{prop:one_block_one_cycle}
  Suppose the partially exposed $G'_p$ contains two
  vertex-disjoint rotating cycles $J_1$ and $J_2$ whose
  vertices are contained in the same block $B\in \mc B$. Then after
  we expose the remaining untested edges of $G'_p$, \aas we can find a
  cycle of length at least $k+1$ in $G'_p[V(B)]$.

\end{proposition}
\begin{proof}
Here, we again use \thmref{thm:menger} to merge cycles. Let $P$ be the
path in $T$ connecting $J_1$ to $J_2$. This path exists and is unique
because $T \cap B$ is a tree. Let $w_1$ be the endpoint of $P$ in $J_1$
and let $w_2$ be the other endpoint of $P$ in $J_2$. We may assume,
without loss of generality, that $w_1$ is the smallest vertex with respect
to $\le_T$ in $J_1$. To see why this assumption can be made, observe that
(P2$\star$) implies that both $V(J_1)$ and $V(J_2)$ induce paths in $T\cap B$.

Since $B$ is a $2$-connected graph, we can use the edges of $B$
to obtain two vertex disjoint paths $P_1,P_2$ connecting
$J_1$ to $J_2$. Let $u_1,v_1$ be the endpoints of $P_1,P_2$
in $J_1$, respectively. Similarly, let $u_2, v_2$ be the
endpoints of $P_1,P_2$ in $J_2$. If neither $u_1$ nor $v_1$
is the pivot of $J_1$, then we can obtain the long cycle
in the following way. By using \propref{prop:rotating_cycles},
we \aas obtain a path $P_3$ of length at least $(2-10\eps)k/3$
between $u_1$ and $v_1$ in $G'_p[V(J_1)]$, and clearly there
exists a path $P_4$ of length at least $(1-4\eps)k/2$ between
$u_2$ and $v_2$ in $J_2$ (just take the longest of the two
paths connecting $u_2$ to $v_2$ in the cycle $J_2$). Putting
together $P_1$, $P_4$, $P_2$, and $P_3$ in that order,
we obtain a cycle of length at least $(7-32\eps)k/6 > k+1$,
thereby proving the proposition.

Otherwise, assume without loss of generality that $u_1$ is
the pivot of $J_1$. One of the edges in the cycle $J_1$ connects
$u_1$ to $w_1$ (recall that $u_1$ is the largest vertex with
respect to $\le_T$, while $w_1$ is the smallest).
The idea now is to modify one of the paths $P_1$ or $P_2$
so that either the endpoint of $P_1$ in $J_1$ is no longer $u_1$,
or the endpoint of $P_2$ is no longer $v_1$, but $w_1$ instead.
To do this, we follow the path $P$ from $w_1$ to $w_2$,
until it hits $P_1$, $P_2$, or $J_2$. If $P$ hits $P_1$ first,
we replace the initial segment of $P_1$ with the initial segment of $P$
as illustrated in \figref{fig:merging_disjoint_cycles}.
If $P$ hits $P_2$ first, we modify $P_2$ similarly. Otherwise,
$P$ never hits $P_1$ or $P_2$, so we can just replace the
whole path $P_2$ by $P$.

\begin{figure}[ht!]
\centering
\begin{tikzpicture}
  [scale=0.8, auto=left, every node/.style={circle, draw, fill=black,
   inner sep=0pt, minimum width=4pt}]
  \node[label=above:$u_1$] (nn_u1) at (-1.94, 1.06 cm) {};
  \node[label=below:$v_1$] (nn_v1) at (-1.94, -1.06 cm) {};
  \node[label=left:$w_1$] (nn_w1) at (-1.5, 0 cm) {};
  \node[label=above:$u_2$] (nn_u2) at (1.94, 1.06 cm) {};
  \node[label=below:$v_2$] (nn_v2) at (1.94, -1.06 cm) {};
  \node[label=right:$w_2$] (nn_w2) at (1.5, 0 cm) {};
  \node (nn_int) at (-0.4, 1.06 cm) {};

  \node[draw=none,fill=none] at (-3, 0cm) {$J_1$};
  \node[draw=none,fill=none] at (3, 0cm) {$J_2$};

  \draw[thick] (-3cm,0cm) circle(1.5cm);
  \draw[thick] ( 3cm,0cm) circle(1.5cm);

  \path[every node/.style={font=\sffamily\small}]
    (nn_u1) edge node[above] {$P_1$} (nn_u2)
    (nn_v1) edge node[below] {$P_2$} (nn_v2)
    (nn_w1) edge[decorate,decoration={snake,amplitude=1.5}]
      node[below] {$P$} (nn_int)
    (nn_int) edge[decorate,decoration={snake,amplitude=1.5}] (nn_w2);

\end{tikzpicture}
\caption{Merging disjoint cycles.}
\label{fig:merging_disjoint_cycles}
\end{figure}
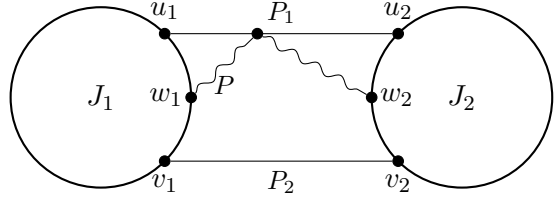

If $P_1$ was modified, we can use the ideas described in the
second paragraph of this proof to obtain the long cycle asymptotically
almost surely. Otherwise, if $P_2$ was modified, we can use a path
connecting $u_1$ to $w_1$ in $J_1$ that uses all of its vertices.
This way we obtain a cycle of length at least $3(1-4\eps)k/2 > k+1$,
and we did not need to test any edge for this case.
\end{proof}

\propref{prop:one_block_one_cycle} dealt with the case of vertex-disjoint
rotating cycles. But what if the cycles intersect? The next proposition
shows that even if the intersection is not empty, it is still possible
to merge the rotating cycles, provided that their intersection is not
too large.

\begin{proposition}

  \label{prop:merge_intersecting_cycles}
  Assume the partially exposed $G'_p$ contains two intersecting
  rotating cycles $J_1$ and $J_2$ contained in the same block
  $B\in \mc B$ whose intersection $J_1 \cap J_2$ has at most $(1-15\eps)k$
  vertices. Then after we expose the remaining untested edges of $G'_p$,
  \aas we can find a cycle of length at least $k+1$ in
  $G'_p[V(B)]$.

\end{proposition}
\begin{proof}
Let $u_1$ and $u_2$ be the pivots of $J_1$ and $J_2$, respectively.
Since $|J_1\cap J_2| \le (1-15\eps)k$, we must necessarily have $u_1\ne
u_2$, as shown in \figref{fig:merging_inter_cycles}.
In fact, if $u_1 = u_2$ then either $V(J_1) \subseteq V(J_2)$ or
$V(J_2) \subseteq V(J_1)$, and hence we would have $|J_1\cap J_2| =
\min\{|J_1|, |J_2|\} > (1-15\eps)k$, which is a contradiction.
Moreover, let $v_1$ and $v_2$ be the smallest vertices in $J_1$ and
$J_2$ respectively, with respect to the order $\le_T$. Furthermore,
let $w$ be the largest vertex in $J_1\cap J_2$ with respect to
the same order. We must have either $v_1 \le_T v_2$ or $v_2 \le_T
v_1$, because otherwise (P2$\star$) would imply that $J_1$ and $J_2$
are disjoint. Assume $v_1\le_T v_2$.
The
intersection $J_1\cap J_2$ comprises the path in $T$ joining $v_2$ to
$w$. We divide the remainder of the proof into two cases.

In the first case we have $|J_1 \cap J_2| < 100 \eps k$.
We can obtain a long cycle $P$ as follows:
we start at $v_1$, traverse the edge to $u_1$, walk
the path in $J_1$ from $u_1$ to $w$ (we choose the path that does
not contain $v_1$), then walk the path in $J_2$ from $w$ to
$u_2$ (again choosing the path the does not contain $v_2$),
move to $v_2$ using an edge from $J_2$, and finish
the cycle with the path from $v_2$ to $v_1$ in $J_1$.
The length of $P$ is at least $|J_1| + |J_2| - 2|J_1\cap J_2| >
k+1$, and we are done.

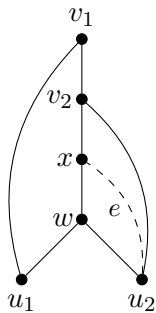
\begin{figure}[ht!]
\centering
\begin{tikzpicture}
  [scale=0.8, auto=left, every node/.style={circle, draw, fill=black,
   inner sep=0pt, minimum width=4pt}]
  \node[label=above:$v_1$] (n_v1) at (0, 3 cm) {};
  \node[label=left:$v_2$] (n_v2) at (0, 2 cm) {};
  \node[label=left:$w$] (n_w) at (0, 0 cm) {};
  \node[label=below:$u_1$] (n_u1) at (-1, -1 cm) {};
  \node[label=below:$u_2$] (n_u2) at ( 1, -1 cm) {};
  \node[label=left:$x$] (n_x) at (0, 1cm) {};

  \path[every node/.style={font=\sffamily\small}]
    (n_v1) edge [bend right] (n_u1)
    (n_v1) edge (n_v2)
    (n_v2) edge (n_x)
    (n_x) edge (n_w)
    (n_w) edge (n_u1)
    (n_v2) edge [bend left] (n_u2)
    (n_w) edge (n_u2)
    (n_u2) edge [bend right, dashed] node[left] {$e$} (n_x);

\end{tikzpicture}
\caption{Merging intersecting cycles.}
\label{fig:merging_inter_cycles}
\end{figure}

In the second case, we have $|J_1 \cap J_2| \ge 100 \eps k$.
Let $X$ be the set of the $\eps k$ vertices $x$ in the path $J_1\cap J_2$
which are closest to $w$ such that $xu_2$ is an untested edge. Such set
$X$ exists because of (P3$\star$) and $|J_1 \cap J_2| \ge 100 \eps k$.
Observe that no vertex in $X$ is more than $4\eps k+\eps k$ vertices
away from $w$, as (P3$\star$) implies.
Next, we expose the untested edges joining $u_2$ to a vertex
in $X$. Asymptotically almost surely we can find a successfully
exposed edge $e=xu_2$. We can now obtain a long cycle $P$ in a way
very similar to what we did before: we start at $v_1$, traverse the
edge to $u_1$, walk the path in $J_1$ from $u_1$ to $w$,
then walk the path in $J_2$ from $w$ to $u_2$,
move to $x$ using the edge $e$ that we recently exposed,
and finish the cycle with the path from $x$ to $v_1$ in $J_1$.
The length of $P$ is at least $|J_1| + |J_2| - |J_1\cap J_2|- 5\eps k >
k+1$, and we are done.
\end{proof}

So far, we only have analyzed the cases where the block has either two
pseudo-cliques or two rotating cycles. To conclude this subsection,
we state a proposition that handles the case when these two different
structures are mixed together in the same block.

\begin{proposition}

  \label{prop:merge_mixed}
  Let $J$ be a rotating cycle of the partially exposed $G'_p$, and let
  $B\in \mc B$ be the unique block containing $J$. Assume that $B$
  contains a pseudo-clique $C\in\mc C$, and that the intersection
  $V(J) \cap C$ has at most $(1-30\eps)k$ vertices.
  Then after we expose the remaining untested edges of $G'_p$, \aas
  we can find a cycle of length at least $k+1$ in $G'_p[V(B)]$.

\end{proposition}
\begin{proof}
The proof is identical to the previous proof and is therefore
omitted.
\end{proof}

\subsection{Step 5: double-counting the poor and the full}
\label{sec:step5}
In this subsection, we study the rotating cycles of $G'_p$. For that purpose,
we assume that the statements of: \propref{prop:small_clique},
\lemref{lem:edgesoutside}, \corref{cor:tested_dfs}, and of
\corref{cor:tested_block} hold. We further assume that \eqref{eqn:size_z1}
holds, and all the edges from $Q_2''$ were tested (as otherwise we would have
a long cycle already). Thus, the reader should bear in mind that any
probabilistic statement in this subsection should be conditioned on the event
that all these assumptions hold.

Let $U$ be the set of all untested edges from $G'_p$ so far. More precisely,
let $U= E(G')\setminus (Q_1 \cup Q_2 \cup Q_3)$. In the next few paragraphs,
we adopt some definitions motivated by the work of Riordan \cite{Riordan}.
We say that a vertex $v$ in a block $B \in \mc B$ is \emph{poor} in $B$ if
the number of descendants (with respect to $T$) of $v$ in $B$ is at most
$\eps k$. Otherwise we say that $v$ is \emph{rich} in $B$. Observe that every
rich vertex in a block $B$ has at least $\eps k$ poor descendants in $B$.
Finally, we say that a vertex $v\in V(B)$ is \emph{full}
in $B$ if the number of vertices $u\in V(B)$ such that $uv\in U$
is at least $(1-\eps)k$.

\obsref{obs:pseudo_path} stated that the vertices of a pseudo-clique
induce a path in the rooted forest $T$. A consequence of this fact
is the following observation.

\begin{observation}

  \label{obs:contribution_poor}
  Let $C\in \mc C$ be a pseudo-clique satisfying $C\not\subseteq W$,
  and let $B\in \mc B$ be unique block containing $C$. The total number
  of vertices in $C$ which are poor in $B$ is at most $\eps k$.

\end{observation}

The next proposition shows how to obtain a rotating cycle from full vertices.

\begin{proposition}

  \label{prop:full_is_rotating}
  Let $v$ be a full vertex in $B\in \mc B$ such that
  the number of edges $vu\in U$ for which $u\in V(B)$ is a descendant
  of $v$ with respect to $T \cap B$ is at most $\eps k$.
  Then by testing some of the untested edges incident
  to $v$, \aas we can obtain a rotating cycle with pivot $v$.

\end{proposition}
\begin{proof}
We would like to remind the reader that we are using the fact
that all the edges of $Q_2''$ were tested (the edges of $Q_2''$
connect vertices at distance greater than $k$ with respect to $T$),
as assumed in the beginning of this subsection. Let $X$ be the set of all
vertices $u\in V(B)$, such that $vu\in U$. Because of the
property of DFS forests stated in \propref{prop:dfs_tree},
we know that for each $u\in X$, $u$ is either a descendant
or an ancestor of $v$. By the hypothesis of the proposition,
the set $X$ has at least $(1-2\eps)k$ ancestors of $v$.
We also know that none of the vertices in $X$ have distance
more than $k$ to $v$ with respect to the tree $T \cap B$.
Let $Y\subseteq X$ consists of the $\eps k$ vertices in $X$
which are ancestors and are as far from $v$ as possible,
with respect to the same distance on the tree $T\cap B$.
Asymptotically almost surely, if we test the edges of $U$
connecting $v$ to vertices in $Y$, we obtain a successfully
tested edge $e = vu$. We claim that the path from $v$ to $u$
in the tree $T\cap B$ together with the edge $e$ forms a
rotating cycle $J$ with pivot $v$. This assertion is clear,
as one can immediately verify that properties (P1$\star$),
(P2$\star$), and (P3$\star$) hold.
\end{proof}

The careful reader will notice that the conditions of
\propref{prop:full_is_rotating} are trivially satisfied when $v$ is full
and poor, hence we have the following corollary.

\begin{corollary}

  \label{cor:fullpoor_is_rotating}
  Let $v$ be a full poor vertex in $B\in \mc B$.
  Then by testing some of the untested edges incident
  to $v$, \aas we can obtain a rotating cycle with pivot $v$.

\end{corollary}

We turn to identify the set of full vertices in the blocks of $\mc B$.
Let $Z_2$ be the set of all vertices $v$ from $G' \setminus Z_1$ such
that $v$ is incident to at least $\frac{\eps k}{3}$ tested edges from
$Q_2 \cup Q_3$. From our assumptions at the beginning of this subsection,
more specifically from Corollaries \ref{cor:tested_dfs} and
\ref{cor:tested_block}, we have
\begin{equation}
  \label{eqn:size_z2}
  |Z_2| \le \frac{15}{\eps p k}\cdot \left(\ell + \frac{n}{k}\right)+12.
\end{equation}

To avoid future issues with double-counting arguments, we would like
to identify the set of vertices $v \in V(G')\setminus (Z_1 \cup Z_2)$, such that
there exists a unique block $B\in \mc B$ for which all but at most
$\eps k$ neighbors of $v$ in $G'$ belong to $V(B)$.
If $v$ is not a cut-vertex of $G'_p$, this is trivial (recall that the
cut-vertices of $G'_p$ are precisely the cut-vertices of $T+(E_1\cup E_3)$).
Otherwise, let $Z_3$ be the set of cut-vertices $v$ from $G'_p$ not in
$Z_1\cup Z_2$ such that $v$ is the smallest (with respect to the order
$\le_T$) of a block containing a pseudo-clique from $\mc C$.
Moreover, let $Z_4$ be the set of all cut-vertices $v$ from $G'_p$
not in $Z_1\cup Z_2 \cup Z_3$, such that there are at least
$\frac{\eps k}{3}$ edges $vw$ in $U$ for which $v$ is the smallest
vertex in the unique block that contains both $v$ and $w$.
We claim the following.

\begin{proposition}

  \label{prop:z3_z4_small}
  $|Z_3| \le |\mc C|$ and $|Z_4| \le \frac{3(\ell+|\mc C|)}{\eps k}$.

\end{proposition}
\begin{proof}
To prove $|Z_3| \le |\mc C|$ note that for each pseudo-clique
$C\in \mc C$ there exists a unique block $B\in\mc B$ such that
$B$ contains $C$. Moreover, there is a unique vertex $v$ which
is the smallest vertex of $B$ with respect to the order $\le_T$.
The map given by $C\mapsto v$ covers every vertex from $Z_3$,
hence $|Z_3| \le |\mc C|$.

To prove the other inequality, observe that
if $vw\in U$, where $v\in Z_4$ and $w$ is a vertex that
belongs to a block $B$ where $v$ is the smallest vertex, then either $w$ is an
outcast vertex, or there exists $C\in \mc C$ such that $w \in C$.
But if $w$ belongs to the pseudo-clique $C$, we claim that $w$ must be the
smallest vertex in the unique block that contains $C$. To see this, first
observe that since $v\not\in Z_3$,  $B$ does not contain $C$. Let $B'\ne B$
be the unique block containing $C$. By \obsref{obs:block_properties},
$w \in V(B) \cap V(B')$ must be the smallest vertex of either $B$ or $B'$.
But because $v$ is the smallest vertex from $B$, we infer that $w$ is the
smallest vertex from $B'$, proving our claim.
Since $w$ is either outcast or the smallest vertex of a block that contains a
pseudo-clique, we must conclude that there are at most $\ell + |C|$ different
choices for $w$.

We claim that for every vertex $w\in V(G')$, there is at most one
edge in $U$ connecting $w$ to an ancestor of $w$ which is the smallest
vertex of some block. To prove this claim, suppose towards contradiction
that there exist two such edges $wx_1$ and $wx_2$. Let $B_1$ and $B_2$
be the corresponding blocks containing $wx_1$ and $wx_2$, respectively.
The intersection of the blocks $B_1$ and $B_2$ contains $w$, but $w$ is
not the smallest in neither of them, contradicting
\obsref{obs:block_properties}, and proving our second claim.

Therefore, there are at most $\ell+|C|$ edges $vw \in U$
such that $v\in Z_4$ and $v$ is the smallest vertex in the unique block
containing both $v$ and $w$. This immediately implies that $|Z_4|
\le \frac{3(\ell+|C|)}{\eps k}$, concluding the proof of the proposition.
\end{proof}

Let $Z= Z_1\cup Z_2 \cup Z_3 \cup Z_4$. Combining \eqref{eqn:size_z1},
\eqref{eqn:size_z2} and \propref{prop:z3_z4_small},
we obtain that \aas
\begin{equation}
  \label{eqn:z_small}
  |Z| \le \frac{1.05n}{k} + \frac{16}{\eps p k}\cdot
  \left(\ell +\frac{n}{k}\right) + 12.
\end{equation}
\begin{proposition}

  \label{prop:not_z_full}
  Let $v \in V(G')$ be an outcast vertex such that $v\not\in Z$.
  Then there exists a unique block $B\in \mc B$ such that
  $v$ is full in $B$. Moreover, for any other block $B'\in \mc B$
  such that $v\in V(B')$ and $B'\ne B$, $v$ is necessarily the smallest
  vertex in $B'$ with respect to the order $\le_T$.

\end{proposition}
\begin{proof}
Since $v$ is an outcast vertex and $v\not\in Z_1$, we know that $\deg_{G'}(v)
\ge \lt(1 -\frac{\eps}{3}\rt) \deg_{G}(v) \ge \lt(1 -\frac{\eps}{3}\rt)k$.
Moreover, because $v\not\in Z_2$, at most $\eps k/3$ edges from $Q_2\cup Q_3$
are incident to $v$. Hence at least $\lt(1-\frac{2\eps}{3}\rt)k$ edges from
$U$ are incident to $v$. Now we split the analysis into two cases:

In the first case, $v$ is not a cut-vertex from $G'_p$. Then there exists a
unique block $B \in \mc B$ such that $v \in V(B)$. Clearly all the edges
in $U$ incident to $v$ are of the form $vu$, for some $u\in V(B)$.
Thus $v$ is full in $B$, and $v$ does not belong to any other block,
concluding the analysis in this case.

In the last case, $v$ is a cut-vertex from $G'_p$. We claim that there exists
a unique block $B \in \mc B$ such that $v$ is not the smallest element in $B$
with respect to the order $\le_T$. To see this, observe that if $v$ is the
smallest vertex in every block in which it belongs, then $v$ must be a
root of the rooted forest $T$, hence $v\in Z_3\cup Z_4$, as $v$ is incident
to more than $\eps k / 3$ edges from $U$. But this is not the case,
therefore there exists at least one block $B$ such that $v\in V(B)$ and $v$
is not the smallest vertex in $B$. By \obsref{obs:block_properties} we
know that such $B$ must be unique. Thus, for all edges in $U$ of the form
$vu$, where $u\not\in V(B)$, the vertex $v$ must necessarily be the smallest
in the unique block that contains both $u$ and $v$. But because
$v\not\in Z_3\cup Z_4$, there can be at most $\eps k /3$ of such edges,
therefore $v$ is full in $B$, concluding the proof of the proposition.
\end{proof}

One immediate consequence of \propref{prop:not_z_full} is the following
corollary.

\begin{corollary}

  \label{cor:not_z_full}
  Let $v\in V(B) \setminus Z$ be a vertex which is not the smallest in
  $B \in \mc B$ with respect to $\le_T$. Then either $v$ is full in $B$,
  or there exists a pseudo-clique $C$ such that $v\in C$ and $B$
  contains $C$.

\end{corollary}
\begin{proof}

If $v$ is an outcast vertex, then since $v$ is not the smallest of
$B$ and $v\not\in Z$, \propref{prop:not_z_full} implies that $v$
must be full in $B$. Otherwise there exists a pseudo-clique $C\in \mc C$
such that $v\in C$. We  claim that $B$ contains $C$. Suppose towards
contradiction that $B$ does not contain $C$. Let $B'$ be the unique
block containing $C$. By \obsref{obs:block_properties}, $v$ must
be the smallest vertex of $B'$, since $v$ is not the smallest
vertex in $B$. Therefore $v$ is the smallest vertex of the block $B'$
which contains the pseudo-clique $C$, hence $v\in Z_3\subseteq Z$,
a contradiction, concluding the proof of the corollary.
\end{proof}

Finally, we turn to the analysis of the poor vertices in the blocks of
$\mc B$. We say that a block $B\in \mc B$ is \emph{good} if the proportion
of vertices from $Z$ in $B$ is at most $\eps/10^6$. We have the following.

\begin{lemma}

  \label{lem:few_poor}
  If the proportion of poor vertices inside a good block $B$ is
  at most $\eps/1000$ then after we expose the remaining untested edges
  of the partially exposed $G'_p$, \aas $G'_p[V(B)]$ contains
  a cycle of length at least $k+1$.

\end{lemma}

Before we prove \lemref{lem:few_poor} we need to prove some
auxiliary results. For that purpose, let us introduce additional
notation.

For each vertex $v\in V(B)$, let $D(v)$ denote the set
of all vertices $u\in V(B)$ which are descendants of $v$ (recall
that every $v$ is a descendant of itself)
and $A(v)$ denote the set of all vertices $u\in V(B)$ which
are ancestors of $v$ with respect to the tree $T\cap B$.
The block $B$ should be clear from the context whenever
we use the notation for ancestors and descendants. We shall
add the subscript ``$\le d$'' to either $D(v)$ or $A(v)$,
such as in the expression $D_{\le d}(v)$, to refer to
the subset obtained by keeping the vertices
at distance at most $d$ from $v$ with respect to the same tree
$T\cap B$. Similarly, we add the superscript ``(p)''/``(r)''
to select only the poor/rich vertices of the indicated set
in the notation, such as in the expression $D^{(p)}(v)$.

We say that a vertex $v\in V(B)$ is \emph{branching} if there exist
at least two distinct rich vertices $u_1, u_2\in V(B)$ which are immediate
descendants of $v$ with respect to $T$.
Similar to what we did previously, we reserve the superscript ``(b)''
to denote the branching vertices of the set under consideration.
We claim the following.

\begin{proposition}

  \label{prop:num_branching}
  For each $v$ such that $D^{(b)}(v)\ne \emptyset$, we have
  $\lt|D^{(p)}(v)\rt| \ge \eps k\lt(\lt|D^{(b)}(v)\rt| + 1\rt)$.

\end{proposition}
\begin{proof}
Assume $D(v)$ contains at least one branching vertex. In this
case $v$ must be rich.
Let $T'$ be the subtree of $T\cap B$ containing all the rich descendants
of $v$ (including itself). Since every branching vertex in $D^{(b)}(v)$
has degree at least $3$ in $T'$ (except possibly the root $v$),
the number of leaves in $T'$ is at least $\lt|D^{(b)})\rt| + 1$.
But every leaf of $T'$ contains at least $\eps k$ poor
descendants in $T\cap B$, thereby proving the proposition.
\end{proof}

\propref{prop:num_branching} yields an upper bound on the total number
of branching vertices in a block, namely it is at most $\frac{1}{\eps k}$
times the number of poor vertices in the same block.

The next proposition allow us to find a structure that resembles
a path with a small number of ``pendant'' vertices in a block with
very few poor vertices.

\begin{proposition}

  \label{prop:good_path}
  If the proportion of poor vertices inside a good block $B$ is
  at most $\eps/1000$ then there exist two vertices $u,v \in V(B)$
  such that $u$ is a descendant of $v$ at distance $30k$ with
  respect to the tree $T\cap B$, and
  \begin{enumerate} \smalllist
    \item the number of vertices in $D(v)\setminus D(u)$ is at
    most $(30 + \eps/10)k$, and
    \item the number of vertices in $\lt(D(v)\setminus D(u)\rt)\cap Z$
    is at most $\eps k/10$.
  \end{enumerate}

\end{proposition}
\begin{proof}
Let $q(x)=\frac{1}{\eps}\cdot\lt|D_{\le d}^{(p)}(x)\rt|+k\cdot
\lt|D_{\le d}^{(b)}(x)\rt|+\frac{1}{\eps}\cdot\lt|D_{\le d}(x)\cap Z\rt|$,
where $d:=40k$. We would like to estimate $q:=\sum\limits_{x~\text{rich}} q(x)$.
We have
\begin{align*}
q &= \sum_{x~\text{rich}}\lt(\sum_{y\in D_{\le d}^{(p)}(x)} 1/\eps
+\sum_{y\in D_{\le d}^{(b)}(x)} k+\sum_{y\in Z\cap D_{\le d}(x)} 1/\eps\rt)\\
&= \sum_{y~\text{poor}} \frac{1}{\eps}\cdot \lt|A_{\le d}^{(r)}(y)\rt|+
   \sum_{y~\text{branching}} k \cdot \lt|A_{\le d}^{(r)}(y)\rt|+
   \sum_{y\in Z} \frac{1}{\eps}\cdot \lt|A_{\le d}^{(r)}(y)\rt| \\
&\le \sum_{y~\text{poor}} \frac{d}{\eps}+
  \sum_{y~\text{branching}} k \cdot d+
  \sum_{y\in Z} \frac{d}{\eps} \\
&\le \frac{d}{\eps}\cdot
   \frac{\eps}{1000}\cdot|B|+ k\cdot d \cdot \frac{1}{\eps k}
   \cdot \frac{\eps}{1000}\cdot|B|+ \frac{d}{\eps}\cdot\frac{\eps}{10^6}
   \cdot |B|< \frac{d}{400}\cdot |B|,
\end{align*}
where for the second-last inequality we used \propref{prop:num_branching}
to estimate the number of branching vertices.
Note that all sums are taken over vertices in $B$.
By averaging, there exists a vertex $v\in V(B)$ such that
$q(v) < d/400=k/10$. In particular, we must have $D_{\le d}^{(b)}(v)
= \emptyset$, $\lt|D_{\le d}^{(p)}(v)\rt| < \eps k/10$ and
$\lt|D_{\le d}(v) \cap Z\rt| < \eps k/10$. Let $d'=30k$.
We claim that for each rich vertex $x\in D_{\le d'}(v)$,
there exists exactly one rich vertex $x'$ which
is an immediate descendant of $x$ with respect to $T\cap B$.
Clearly there are no two such vertices $x'$, since otherwise
$x$ would be branching, and this cannot happen because
$D_{\le d}^{(b)}(v)=\emptyset$. To finish the proof of the claim,
notice that if all the immediate descendants of $x$ were poor,
then $D(x) = \{x\} \cup D^{(p)}(x) = \{x\} \cup D^{(p)}_{\le \eps k}(x)$,
which together with $d > d' + \eps k$ implies that
$\lt|D_{\le d}^{(p)}(v)\rt| \ge\lt|D_{\le \eps k}^{(p)}(x)\rt| \ge \eps k$,
a contradiction.

By the claim we proved in the previous paragraph, we know that
the set $D_{\le d'}^{(r)}(v)$ induces a path in $T$.
Let $u$ be the unique rich vertex in $D(v)$ at distance exactly
$d'$. We claim that the pair $u,v$ satisfies the conditions
stated in the proposition. For the first condition, observe that
$D(v)\setminus D(u) \subseteq D_{\le d'}^{(r)}(v) \cup D_{\le d}^{(p)}(v)$,
hence clearly $|D(v)\setminus D(u)| < d' + \eps k / 10$.
For the second condition, we have that $Z \cap \lt(D(v)\setminus D(u)\rt)
\subseteq D_{\le d}(v)\cap Z$, hence $\lt|Z \cap \lt(D(v)\setminus D(u)\rt)\rt|<
\eps k / 10$, finishing the proof of the proposition.
\end{proof}
We have the necessary tools to prove \lemref{lem:few_poor}.
\begin{proof}[Proof of \lemref{lem:few_poor}{}]
We assume, without loss of generality, that $B$ contains at most one
pseudo-clique from $\mc C$, as otherwise \propref{%
prop:one_block_one_clique} would already imply the conclusion
of this lemma.

We start the proof by applying \propref{prop:good_path} to
$B$, thus obtaining the pair $u,v$. Let $P$ be the path
between $u$ and $v$ in $T$. We have $|P| = 30 k$, the number
of vertices in $V(P)\cap Z$ is at most $\eps k /10$, and the
number of ``pendant'' vertices from $P$ is at most $\eps k /10$.
In particular, for each vertex $w \in V(P)$ at distance
at least $k$ from $u$ with respect to $P$, there are
at most $\eps k/10$ edges in $U$ from $w$ to one of its
descendants not in $P$.

We redefine $P$ to be the subpath of length $28k$ obtained by
removing the two segments of length $k$ closest to the
two endpoints from the original path. For each vertex $x\in
V(P)\setminus Z$, we know (by \corref{cor:not_z_full}) that
either $x$ is full in $B$, or $B$ contains a pseudo-clique
$C\in \mc C$ such that $x\in C$. But we know, by our initial
assumption, that there is at most one such $C$, and if it
exists, then $V(P) \cap C$ should be a segment of
$P$ (because pseudo-cliques induce paths in $T$,
see \obsref{obs:pseudo_path}). Hence, we can always find in $P$
two disjoint segments $L_1, L_2$, each of length $8k$, such that
for every $x \in \lt(V(L_1)\cup V(L_2)\rt) \setminus Z$, $x$ is
full in $B$ and the distance between these segments along the path $P$ is at least $k$.
In other words, almost all vertices from $L_1\cup L_2$
are full in $B$.

Let us divide the rest of the proof into two cases. In the first case,
we assume that there exist two vertices $x_1\in L_1$ and $x_2\in L_2$,
both full in $B$, such that for each $i\in \{1,2\}$, there are at most
$\eps k$ descendants $w$ of $x_i$ in $B$ for which $wx_i$ is an
edge in $U$. By \propref{prop:full_is_rotating}, then \aas we can
obtain two rotating cycles, with pivots $x_1$ and $x_2$ and by
\propref{prop:one_block_one_cycle} we can merge these two
disjoint rotating cycles and obtain the desired long cycle,
proving the lemma in this first case.

In the second case, we assume that there is no such pair of
vertices $x_1, x_2$. Hence we might also assume that, without loss of generality,
for each full vertex $x$ in $L_1$, there exist at least $\eps k$
descendants $y$ of $x$ for which $xy \in U$. Out of these descendants,
at most $\eps k / 10$ do not belong to $P$ (recall that the number of
``pendant'' vertices from $P$ is at most $\eps k / 10$). Thus $x$
sends at least $9\eps k / 10$ untested edges to its descendants in $P$.
Furthermore, at most $4\eps k / 5$ of these descendants are of
distance at most $4 \eps k /5$, thus there are at least $\eps k / 10$
edges in $U$ connecting $x$ to one of its descendants in $P$ at distance
at least $4\eps k / 5$ from $x$.
Observe that we can a.a.s. obtain a cycle of length at least $4\eps k / 5$
by testing the $\eps k / 10$ edges in $U$ going from a full vertex to its
descendants on the path $P$. The key idea in what comes next is to
merge $O(1/\eps)$ of these small cycles.

In $P$, and hence in $L_1$, there can be at most $\eps k/10$
non-full vertices. This implies that for each subsegment $L$ in $L_1$
of length $\eps k/5$, at least $\eps k / 10$ of its vertices are full,
thus there exists a set $E_{L}$ of at least $\eps^2 k^2/100$ edges
in $U$ of the form $xy$, where $x\in V(L)$, and $y \in V(P)$ is a descendant of
$x$ at distance at least $4\eps k /5$. Clearly the distance between $x$
and $y$ in $P$ is at most $k$, as all edges of $Q_2''$ were tested
according to our assumption in the beginning of this subsection.
By the union bound and by \lemref{lem:chernoff}, if we test the
edges in $E_{L}$ for every segment $L$ in $L_1$, \aas we can find
one successfully exposed edge in each $E_{L}$.

To obtain the long cycle is straightforward.
We start with the a segment $L^{(1)}$ of length $\eps k /5$
containing the endpoint of $L_1$
which is smallest with respect to the order $\le_T$. In this segment,
we can find an edge $u_1v_1\in E_{L^{(1)}}$ which was successfully exposed,
where $u_1\in V(L^{(1)})$. We then proceed recursively for each $j$ as
follows: let $L^{(j+1)}$ be the segment of $L_1$ of length $\eps k / 5$
whose smallest vertex with respect to $\le_T$ is the
ancestor of $v_j$ at distance $\eps k / 5 + 1$ (hence $L^{(j+1)}$ does not
contain $v_j$). Then choose an edge
$u_{j+1}v_{j+1}\in E_{L^{(j+1)}}$ which was successfully exposed, where
$u_{j+1}\in V\lt(L^{(j+1)}\rt)$. Repeat this process while $v_j$ has
distance at most $5k$ from the smallest vertex from $L_1$. Assume
the last segment chosen was $L^{(t)}$. The cycle $J$ we seek can be
easily seen from \figref{fig:few_poor_long_cycle}.

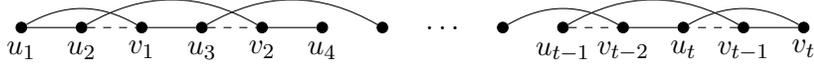
\begin{figure}[ht!]
\centering
\begin{tikzpicture}
  [scale=0.8, auto=left, every node/.style={circle, draw, fill=black,
   inner sep=0pt, minimum width=4pt}]
  \node[label=below:$u_1$] (u1) at (-3, 0 cm) {};
  \node[label=below:$u_2$] (u2) at (-2, 0 cm) {};
  \node[label=below:$v_1$] (v1) at (-1, 0 cm) {};
  \node[label=below:$u_3$] (u3) at (0, 0 cm) {};
  \node[label=below:$v_2$] (v2) at (1, 0 cm) {};
  \node[label=below:$u_4$] (u4) at (2, 0 cm) {};
  \node (v3) at (3, 0cm) {};
  \node[draw=none, fill=none] (empty) at (4, 0cm) {$\ldots$};
  \node (uu1) at (5, 0 cm) {};
  \node[label={[yshift=4pt]below:$u_{t-1}$}] (uu2) at (6, 0 cm) {};
  \node[label={[yshift=4pt]below:$v_{t-2}$}] (vv1) at (7, 0 cm) {};
  \node[label={below:$u_t$}] (uu3) at (8, 0 cm) {};
  \node[label={[yshift=4pt]below:$v_{t-1}$}] (vv2) at (9, 0 cm) {};
  \node[label={below:$v_t$}] (uu4) at (10, 0 cm) {};

  \path[every node/.style={font=\sffamily\small}]
    (u1) edge (u2)
    (u1) edge[bend left] (v1)
    (u2) edge[dashed] (v1)
    (u2) edge[bend left] (v2)
    (u3) edge (v1)
    (u3) edge[dashed] (v2)
    (u3) edge[bend left] (v3)
    (v2) edge (u4)
    (uu4) edge (vv2)
    (uu4) edge[bend right] (uu3)
    (vv2) edge[dashed] (uu3)
    (vv2) edge[bend right] (uu2)
    (uu3) edge (vv1)
    (vv1) edge[bend right] (uu1)
    (uu2) edge[dashed] (vv1);

\end{tikzpicture}
\caption{The long cycle $J$ formed by the solid lines.}
\label{fig:few_poor_long_cycle}
\end{figure}
The reader can check that the distance between $u_{j}$
and $v_{j-1}$ is always greater than the distance between
$u_{j}$ and $v_{j-2}$ with respect to the path $P$. In fact,
we always jump ``downwards'' by at least $4\eps k / 5$ and move
``upwards'' by $\eps k / 5 + 1$. So the total length
of $J$ is at least two thirds of the distance between $u_1$ and $v_t$,
therefore $|J| \ge 3k$, which proves the lemma.
\end{proof}

Now that we have studied the case when the number of poor vertices
is small, it is natural to consider the case where this number is
sufficiently large. In the next lemma, we show the existence of the
desired cycle in this situation.

\begin{lemma}

  \label{lem:many_poor}
  If the number of poor vertices not in $Z$ inside a block $B$ is
  at least $200 \eps k$, then after we expose the remaining untested
  edges of $G'_p$, \aas $G'_p[V(B)]$ contains a cycle of length
  at least $k+1$.

\end{lemma}

Before we prove the previous lemma, we list and prove some technical results.
The first one is Theorem 3.1 (ii) proved by Krivelevich, Lee, and
Sudakov~\cite{KrivLeeSud}.

\begin{theorem}
\label{thm:bipartite_long_path}
  Let $G$ be a bipartite graph of minimal degree at least $k$. Then
  asymptotically almost surely $G_p$ contains a path of length
  $(2+o(1))k$ whenever $p = \frac{\omega_k(1)}{k}$, for any
  function $\omega_k(1) < k$ that tends to infinity with $k$.
\end{theorem}

The second technical result is inspired by the same paper \cite{KrivLeeSud}.

\begin{proposition}
  \label{prop:benny}
  Let $B\in \mc B$ be a block. Suppose there exists a pair $(P,X)$,
  where $P$ is a path in $T\cap B$ of size $(1-5\eps)k\le |P|\le (1+\eps)k$,
  and $X \subseteq V(B)\setminus V(P)$ is a set of at least $180\eps k$
  vertices, such that for each $v\in X$ there are at least $(1-80\eps)k$
  vertices $u\in V(P)$ satisfying $uv\in U$. After exposing
  all the edges in $U$ connecting $P$ to $X$, \aas we can find a
  cycle in $G'_p[V(B)]$ of length at least $k+1$.
\end{proposition}

\begin{proof}
The main idea here goes along the lines of the proof of Theorem 1.2 (Case 1)
in the above mentioned paper~\cite{KrivLeeSud}. Let $w$ be an arbitrary vertex in $X$.
We test all the edges in $U$ connecting $w$ to $P$, and \aas we obtain
two successfully tested edges $wu_1$ and $wu_2$, such that $u_1$ and $u_2$ are
at distance at least $(1-82\eps)k$ with respect to $P$. We redefine $P$ as
the segment of itself connecting $u_1$ to $u_2$. We know now that $(1-82\eps)k
\le |P| \le (1+\eps)k$, and that for each vertex $v \in X\setminus \{w\}$,
there exists a set of at least $(1-163\eps)k$ vertices $u\in V(P)$ such that
$uv\in U$.

Let $Y$ be an arbitrary subset of $X\setminus \{w\}$ of size
$175\eps |P| < 180\eps k$. We partition $P$ into $\frac{1}{175\eps}$
segments $P_1,\ldots,P_{1/175\eps}$, each of length $175\eps |P|$.
By an averaging argument, there exists an interval $P_i$
for which $e(P_i, X) \ge (1-163\eps)|P_i||Y|$. Consider a bipartite
graph $\Gamma$ formed by the edges of $U$ with the vertex set being the
union of the two parts $P_i$ and $Y$. Note that the number of non-adjacent
pairs is at most $163\eps |P_i||Y|$ (also note that $|Y|=|P_i|$). We repeatedly
remove vertices of degree at most $(1-40\eps^{1/2})|Y|$ from $\Gamma$. As
long as the total number of deleted vertices is at most $20\eps^{1/2}|Y|$,
each deletion accounts for at least $20\eps^{1/2}|Y|$ non-adjacent pairs
of $\Gamma$. So, if we continue this removal process for at least
$20\eps^{1/2}|Y|$ vertices, the total number of non-adjacent pairs we removed
from $\Gamma$ is at least $400\eps|Y|^2$, which is a contradiction. Thus,
this process must stop before we remove $20\eps^{1/2}|Y|$ vertices,
and we obtain a subgraph $\Gamma_1$ of minimum degree at least
$(1-40\eps^{1/2})|Y|$.

Let $P_{i,0}$ and $P_{i,1}$ be the two segments of $P_i$ of length
$45\eps^{1/2}|Y|$ closest to the two endpoints of $P_i$. Even after
removing the vertices in $P_{i,0} \cup P_{i,1}$ from $\Gamma_1$, we
are left with a graph $\Gamma_2$ of minimum degree at least
$(1-40\eps^{1/2})|Y|-90\eps^{1/2}|Y|>\frac{9}{10}|Y|$.

By \thmref{thm:bipartite_long_path}, $(\Gamma_2)_p$ \aas contains
a path of length at least $2(\frac{9}{10}+o(1))|Y| > \frac{5}{3}|Y|$.
By removing at most two vertices, we may assume that the endpoints
$x$ and $y$ of this path are both in $Y$. Since $\Gamma_1$ has
minimum degree at least $(1-40\eps^{1/2})|P_i|$, both of these
endpoints have at least $5\eps^{1/2}|P_i|> \eps^{3/2}k$
neighbors in both $P_{i,0}$ and $P_{i,1}$. By \lemref{lem:chernoff},
$(\Gamma_1)_p$ \aas contains two edges $xv_0$ and $yv_1$, where
$v_0\in P_{i,0}$ and $v_1\in P_{i,1}$. We found a path in $(\Gamma_1)_p$
of length at least $\frac{5}{3}|Y|$, which starts at $v_0$ and ends
at $v_1$, and uses only vertices from $Y\cup (P_i\setminus (P_{i,0} \cup
P_{i,1}))$ as internal vertices. Together with the path $P$ and the
two edges $wu_1$ and $wu_2$, we obtain a cycle of length at least
$|P|-|P_i|+\frac{5}{3}|Y| > k + 1$, concluding the proof of the
proposition.
\end{proof}
Our last auxiliary statement studies the set of poor vertices
inside a block.
\begin{proposition}

  \label{prop:poor_in_block}
  Let $B$ be a block. Suppose $x_1$ and $x_2$ are two full poor vertices
  in $B$. If $x_1$ and $x_2$ have distance at least $30\eps k$ with respect
  to $T$, then after we expose the remaining untested edges of $G'_p$,
  \aas $G'_p[V(B)]$ contains a cycle of length at least $k+1$.

\end{proposition}
\begin{proof}
\propref{prop:full_is_rotating} guarantees that, by testing some edges
in $U$ incident to $x_1$ and $x_2$, we \aas will obtain two rotating
cycles $J_1$ and $J_2$ with pivots $x_1$ and $x_2$, respectively. If $J_1$
and $J_2$ are disjoint, by \propref{prop:one_block_one_cycle} we
can merge them, and obtain the desired long cycle. Otherwise,
if $J_1$ and $J_2$ intersect, then the unique path
in $T$ from $x_1$ to $x_2$ is contained in $J_1\cup J_2$ and is
edge-disjoint from $J_1\cap J_2$. Thus
\begin{equation}
  \label{eqn:size_inter}
  |J_1| + |J_2| - |J_1\cap J_2| = |J_1\cup J_2| \ge
  |J_1\cap J_2| + 30\eps k,
\end{equation}
which implies $|J_1\cap J_2| <(1-15\eps) k$.
By \propref{prop:merge_intersecting_cycles}, we can merge these
two intersecting rotating cycles, thereby proving the
proposition.
\end{proof}

We turn to prove \lemref{lem:many_poor}.

\begin{proof}[Proof of \lemref{lem:many_poor}{}]
\obsref{obs:contribution_poor} says that each pseudo-clique can contribute
at most $\eps k$ poor vertices to the block in which it is contained.
Moreover, if $B$ contains more than one pseudo-clique,
\propref{prop:one_block_one_clique} would already imply the conclusion
of this lemma. Thus, we may assume that among the poor vertices not
in $Z$, at least $199\eps k$ of them do not belong to pseudo-cliques
that are contained in $B$. Let $X$ denote the set such vertices.
We have $|X|\ge 199\eps k$, and clearly $X$ does not contain the
smallest vertex of $B$ with respect to $\le_T$ (otherwise the
smallest element would be poor, and hence $|B| \le \eps k$).

By \corref{cor:not_z_full}, each vertex in $X$ must necessarily
be full in $B$. If there exist two full poor vertices $v_1$ and $v_2$
at distance at least $30\eps k$ with respect to $T$, we can obtain
the long cycle by using \propref{prop:poor_in_block}. Thus
we may assume that all the poor full vertices of $B$ are close to each
other, i.e., have distance at most $30\eps k$. Let $v$ be any full poor
vertex, and let $P$ be the path containing all the ancestors of $v$ in
$B$ at distance less than $k$. Clearly $(1-2\eps)k \le |P| \le k$.
For every vertex $v'\in X$, since $v$ and $v'$ are at distance at most
$30\eps k$, there exists at least $(1-80\eps)k$ pairs $uv'\in U$, where
$u\in P$ (recall that no edge in $U$ connects pairs of vertices at
distance larger than $k$). Thus we can apply \propref{prop:benny}, and
\aas obtain a cycle of length at least $k+1$ in $G'_p[V(B)]$, finishing
the proof of the lemma.
\end{proof}

The last case to be solved is when the number of poor vertices
is not too large and not too small. The next lemma investigates this case.
\begin{lemma}

  \label{lem:some_poor}
  Suppose the proportion of poor vertices inside a good block $B\in \mc B$ is
  at least $\eps/1000$, but the total number of poor vertices in $B\setminus Z$
  is at most $200\eps k$. We have either
  \begin{itemize} \smalllist
  \item $B$ contains a pseudo-clique $C$ and all the other
    vertices in $V(B)\setminus C$ belong to $Z$, or
  \item after testing all the edges of $U$ joining two vertices from
    $B$, \aas $G'_p[V(B)]$ contains a cycle of length at least $k+1$.
  \end{itemize}

\end{lemma}

As before, we need some technical statements in preparation for the
proof of \lemref{lem:some_poor}. The first statement strengthens
\propref{prop:poor_in_block}.

\begin{proposition}

  \label{prop:poor_in_block2}
  Let $B$ be a block. Suppose $x_1, x_2\not\in Z$ are two vertices
  in $B$, such that for each $i\in\{1,2\}$, there are at most
  $\eps k$ edges in $U$ connecting $x_i$ to one of its descendants
  in $B$. If $x_1$ and $x_2$ have distance at least $60\eps k$ with respect
  to $T$ and are not comparable with respect to $\le_T$, then after we
  expose the remaining untested edges of $G'_p$, \aas $G'_p[V(B)]$
  contains a cycle of length at least $k+1$. In particular, the
  conclusion of this proposition also holds if $x_1$ and $x_2$ are
  two poor vertices in $B$ which are not in $Z$.

\end{proposition}
\begin{proof}
First observe that the smallest vertex in $B$ does not belong to
the set $\{x_1, x_2\}$, because otherwise $x_1$ and $x_2$ would be
comparable with respect to $\le_T$. Let $L$ be the union of all
pseudo-cliques $C\in \mc C$ which are contained in $B$.

In the first case, both $x_1$ and $x_2$ belong to $L$.
Let $C_1, C_2\in \mc C$ be the pseudo-cliques contained in $B$
such that $x_1\in C_1$ and $x_2 \in C_2$. We claim that $C_1\ne C_2$.
To prove our claim, assume that $C_1=C_2$. But by \obsref{obs:pseudo_path},
we know that the vertices of $C_1$ induce a path in $T$, and
hence in $T\cap B$. But this would imply that $x_1$ and $x_2$
are comparable with respect to $\le_T$, contradicting the
hypothesis of the proposition. Hence $C_1\ne C_2$. But
by \propref{prop:one_block_one_clique}, we can merge the
two cycles in the pseudo-cliques and obtain the desired
long cycle. This finishes the analysis of the first case.

In the second case, both $x_1$ and $x_2$ do not belong to $L$.
By \corref{cor:not_z_full}, both $x_1$ and $x_2$ are full in $B$,
since they do not belong to $L \cup Z$. After testing some
of the edges of $U$ incident to $x_1$ and $x_2$ we can obtain
two rotating cycles $J_1$ and $J_2$ respectively, as
\propref{prop:full_is_rotating} assures. Either $J_1$ is
disjoint from $J_2$ or their intersection is of size at most
$(1-30\eps)k$ (we use the same strategy as in \eqref{eqn:size_inter}
to estimate the size of $J_1\cap J_2$). In any case, by either
\propref{prop:one_block_one_cycle} or \propref{prop:merge_intersecting_cycles},
\aas we can obtain the long cycle after we test the remaining edges
of $U$ in $B$, concluding the analysis of this case.

In the last remaining case, we assume that $x_1\in L$ but
$x_2\not\in L$. As before, we know that $x_2$ is full in $B$.
Using \propref{prop:full_is_rotating}, we \aas obtain a rotating
cycle $J$ in $B$ for which $x_2$ is its pivot.
Moreover, since $x_1\in L$, we also know that $x_1$ belongs to a
pseudo-clique $C\in \mc C$ which is contained in $B$. We claim
that $|V(J) \cap C| < (1-30\eps)k$. If $V(J)$ and $C$ are disjoint,
then our claim is trivially true. Otherwise, if they intersect,
then since there is a cycle $J'$ in $B$ containing the vertices of
$C\setminus W$, such that $J'\cap T$ is a path (\obsref{obs:pseudo_path}).
By a calculation analogous to \eqref{eqn:size_inter}, we obtain
that $|J\cap J'| = |V(J) \cap C| < (1-30\eps)k$, proving our
claim. Finally, we finish the proof of this proposition with
a final application of \propref{prop:merge_mixed} to obtain
the long cycle.
\end{proof}

Our second auxiliary result allows us to estimate the number
of poor vertices in a block.

\begin{proposition}

  \label{prop:num_poor_distance}
  Let $\delta>1$, $B$ be a block, $v$ be a poor vertex in $B$, and let $P$
  be the unique path from $v$ to the smallest vertex from $B$
  with respect to $\le_T$. If $|V(B)\cap Z| < \eps k$, and every
  poor vertex of $B$ not in $Z\cup \{v\}$ is at distance at most
  $\delta\eps k$ from $v$, then there is no rich vertex in $B$ outside
  of $P$ at distance at least $\delta \eps k$ from $v$. Furthermore,
  $B$ contains at least $\frac{|B|-|P|}{\delta+1}$ poor vertices.

\end{proposition}
\begin{proof}
Let $X$ be the set of poor vertices in $B$. To proof of the first
part of the proposition goes by contradiction. Assume that there exists
a rich vertex $u$ outside of $P$ at distance at least $\delta \eps k$
from $v$. Observe that all vertices from $D(u)\cap X$ must have distance
greater than $\delta \eps k$ from $v$. On the other hand, $|D(u) \cap X|
\ge \eps k$ (since every rich vertex has at least $\eps k$ poor
descendants), and every vertex in $D(u)\cap X$ must belong to $Z$ by the
assumptions of the proposition. We then have a contradiction, because
$|D(u)\cap X| \ge \eps k > |V(B)\cap Z|$. This contradiction
proves the first statement of the proposition.

For each $u \not\in V(P)$ such that $u$ has an immediate ancestor
in $P$, we shall prove that
\begin{equation}
  \label{eqn:poor_proportion}
  |D(u)\cap X| \ge \frac{|D(u)|}{\delta + 1}.
\end{equation}
We identify the set $F$ of the rich vertices in $D(u)$ that have no rich
descendant. For every $w\in F$, since $w$ is rich, we have $|D(w)|> \eps k$,
and all vertices in $D(w)$ are poor, except for $w$ itself. Hence
$|D(u)\cap X| \ge \eps k |F|$. Furthermore, no rich vertex in $D(u)$ is at
distance larger than $\delta \eps k$ from $u$, since it would have distance
at least $\delta \eps k$ from $v$ as well, which is impossible.
In other words, every rich vertex in $D(u)$ belongs to some path in $T$
from a vertex in $F$ to $u$. Since all these paths have length at most
$\delta \eps k$, the total number of rich vertices in $D(u)$
is at most $\delta \eps k |F|$, thereby proving \eqref{eqn:poor_proportion}.
Therefore the total number of poor vertices in $B$ is at least
$\frac{|B| - |P|}{\delta + 1}$, concluding the proof of the proposition.
\end{proof}

The last auxiliary result is to handle the case of a block containing
a long path and few vertices in $Z$. The statement is as follows.

\begin{proposition}

  \label{prop:long_path_good}
  Let $B$ be a block, $v$ be a poor vertex in $V(B)\setminus Z$, and let $P$
  be the unique path from $v$ to the smallest vertex from $B$
  with respect to $\le_T$. If $|V(B)\cap Z| < 3\eps k/10$, $|P|\ge
  (1+1000\eps)k$, and every poor vertex of $B$ not in $Z\cup \{v\}$
  is at distance at most $100 \eps k$ from $v$, then after testing
  the remaining untested edges from $G'_p$, \aas we can find a cycle
  of length at least $k+1$ in $G'_p[V(B)]$.

\end{proposition}
\begin{proof}
Let $L$ be the union of all pseudo-cliques that are contained in $B$.
Observe that $L$ is the union of at most one pseudo-clique, as otherwise
we would have a cycle of length at least $k+1$ by
\propref{prop:one_block_one_clique}.

We claim that either $v\in L$, or after possibly testing few edges
from $U$ incident to $v$, we \aas obtain a rotating cycle $J$ with pivot
$v$. To see this, observe that if $v\not\in L$, then
\corref{cor:not_z_full} implies that $v$ must be full in $B$.
Using \corref{cor:fullpoor_is_rotating}, we \aas obtain
such a rotating cycle $J$, proving our claim. In any case, we can
assume that $v$ belongs to a cycle (either because $v$ belongs to
a pseudo-clique $C$ contained in $B$ or because it is the pivot of
a rotating cycle $J$) of length at least $(1-5\eps)k$.

Suppose there is vertex $w$ in $P$ at distance at
least $100\eps k$ from the endpoint $v$ of $P$ that satisfies:
\begin{enumerate}[(i)] \smalllist
\item $w\not\in Z\cup L$, and $w$ is not the smallest vertex in $B$, and
\item there are at most $\eps k$ edges of $U$ connecting $w$ to
one of its descendants.
\end{enumerate}
By (i) and by \corref{cor:not_z_full}, we know that $w$ must be full in $B$,
and hence by \propref{prop:full_is_rotating}, after testing the far-reaching
edges in $U$ incident to $w$, we \aas obtain a rotating cycle $J'$
with pivot $w$. We could then merge $J'$ with the large cycle containing
$v$ (which could be either from a pseudo-clique $C\in \mc C$ if $v\in L$,
or from the rotating cycle $J$ with pivot $v$, if $v\not\in L$)
by using \propref{prop:merge_mixed}, hence obtaining
the desired long cycle, and we would be done.

From the discussion in the last paragraph, we may assume without
loss of generality that there is no vertex $w$ satisfying both (i) and (ii).
Thus every vertex $w\in V(P)$ at distance at least $100\eps k$ from $v$
satisfying (i) must have at least $\eps k$ descendants $u$ in $B$ such that
$uw\in U$. Out of these descendants, at most $|Z \cap V(B)|< 3\eps k/10$
do not belong to $P$. To see this, observe that $w$ does not have a rich
descendant $u$ outside of $P$, as otherwise it would contradict the conclusion
of \propref{prop:num_poor_distance} (applied with $\delta = 100$).
But every poor descendant of $w$ is at distance at least $100\eps k$
from $v$, hence it must belong to $Z\cap V(B)$. So the total number
of descendants of $w$ outside $P$ is at most $3\eps k/10$. In particular,
$w$ has at least $7\eps k/10$ descendants $u$ such that $u\in V(P)$ and
$uw\in U$.

Note that apart from the vertices in $L$, most vertices in $P$ are at
distance at least $100\eps k$ from $v$ and satisfy (i), because
$|V(B)\cap Z| \le 3\eps  k /10$. But from each vertex $w$ in $P$ satisfying
(i), there are at least $7\eps k/10$ edges in $U$ connecting $w$ to one of its
descendants in $P$. At most $13\eps k/20$ of these edges connect $w$ to a vertex
in $P$ at distance at most $13\eps k/20$ from $w$. Hence there are at least
$\eps k/20$ edges in $U$ connecting $w$ to one of its descendants at
distance at least $13\eps k/20$. If we test these edges, \aas
we can find a successfully tested edge connecting $w$ to one of its
deep descendants in $P$, thus forming a small cycle of length
at least $13\eps k / 20$. We can now use the same technique as
in the proof of \lemref{lem:few_poor} to finish the proof of the
proposition (see \figref{fig:few_poor_long_cycle}).
In the next few paragraphs, we briefly sketch this technique.
We also remark that in our case we only need to merge constantly
many small cycles, which simplifies the union-bound argument.

The idea is to start at a full vertex $w_0$ at distance
between $(1+998\eps)k$ and $(1+999\eps)k$ from $v$, and repeat the
following loop. For each $i=0,1,2,\ldots$, by testing some edges of $U$
incident to $w_i$, we \aas can find a neighbor $w_i'$ of $w_i$ at
distance at least $13\eps k /20$ from $w_i$ which is a descendant
of $w_i$ on $P$. Then we go ``upwards'' (in direction to the smallest
vertex from $B$) the path $P$ starting from $w_i'$ until we find another
full vertex $w_{i+1}$. Recall that we need to go ``upwards'' at most
$3\eps k /10$ vertices to reach this full vertex, as long as we move
entirely outside of $L$. Also observe that the small cycles do not
``double-overlap'', as $13\eps k /20 > 2\cdot (3\eps k /10)$. We repeat
the loop until we either hit the interior of $J$ (if it exists),
or a vertex from $L$ which is not the smallest vertex in $L$.

Recall that $v$ belongs to a cycle, which could be formed by vertices from
either $J$ or $L$. Since this cycle has size between $(1-5\eps)k$ and $(1+\eps)k$,
we must necessarily stop this procedure after constantly many iterations
of the loop. More precisely, if $T$ denotes the time we stopped, then
$T<\frac{999\eps k}{7\eps k/20}<300$.
Moreover, at the very last step, the vertex $w_T'$ either
belongs to the interior of $J$ (if it exists) or is in $L$ (but is not the
smallest vertex in the pseudo-clique). In the first case, we can
close the cycle we are forming with $J$, because $v$ has \aas a
neighbor which is an ancestor of $w_T'$ at distance at
most $5\eps k$ from $w_T'$. In the latter case, when $w_T'\in L$,
from \propref{prop:merge_mixed} one can deduce that the smallest vertex
from $L$ in the block is at distance at most $(1+40\eps)k$ from $v$.
In addition, by \propref{prop:rotating_cliques}, we can close the
cycle using the majority of vertices from $L$. More specifically,
there is a path in $G'_p[L]$ of length at least $(1-20\eps)k$ which
connects $w_T'$ and the smallest vertex from $L$.
Therefore, regardless of what happens in the last iteration of our procedure,
the merged cycle has always size at least $(1+998\eps)k-45\eps k -20\eps k
-300 \cdot \frac{3\eps k}{10}>k$, thereby proving the proposition.
\end{proof}

We are ready to prove \lemref{lem:some_poor}.
\begin{proof}[Proof of \lemref{lem:some_poor}]
Let $X$ be the set of poor vertices of $B$, and let $Y=X\setminus Z$.
Since $|Z\cap V(B)| \le 10^{-6}\eps |B|$ (because $B$ is good) and
$|X| \ge 10^{-3}\eps |B|$, we must conclude that $|Y| \ge \eps |B|/1100$.
We also know that $|Y| \le 200 \eps k$, and this implies that
$|B| < 3\cdot 10^5 k$ (we will improve this bound later),
thus $|Z\cap V(B)| < 3 \eps k / 10$. At last, we have $|X| \le |Y| +
|Z \cap V(B)| \le 201\eps k$.

Fix a vertex $v\in Y$ arbitrarily. Using \propref{prop:poor_in_block2}
one can see that $v$ has distance (with respect to $T$) of at most
$60\eps k$ from any other vertex from $Y$, as otherwise we would
obtain the long cycle and the second conclusion of the lemma would hold.

Let $P$ be the path in $T$ joining $v$ to the smallest vertex $v_0$
of $B$ with respect to $\le_T$. By applying \propref{prop:num_poor_distance}
for $\delta = 60$, we obtain that the total number of poor vertices in
$B$ is at least $(|B| - |P|)/61$, which implies that $|B| \le |P| + 61|X|\le
|P| + 15000\eps k$.

We might assume then that $|P| < (1 + 1000\eps)k$. This is because
the conclusion of the lemma would be true otherwise, as
\propref{prop:long_path_good} shows. In particular, we must have
$|B| < (1 + 16000\eps)k$. If $|V(B)\setminus Z|>(1+\eps)k$,
then we claim that $\Gamma:=G'[V(B)\setminus (Z\cup \{v_0\})]$
is a graph with minimum degree at least $(1-5\eps)k$. Indeed,
every vertex in $\Gamma$ which does not belong to a pseudo-clique
contained in $B$ is full in $B$, and every vertex in $\Gamma$
which does belong to a pseudo-clique contained in $B$ has
degree at least $(1-5\eps)k$ in $G'[V(B)\setminus (Z\cup \{v_0\})]$
because $|Z\cap V(B)|< 3\eps k / 10$ and every pseudo-clique not in
completely inside the waste lost at most $\eps k/2$ vertices to $W$,
and the claim follows.

Thus $\Gamma$ is a graph with minimum degree at least $(1-6\eps)k$
satisfying $(1+\eps)k \le |\Gamma| \le (1+16000\eps)k$. However, as
assumed in the beginning of this subsection, such graph cannot exist,
because it would violate the assumption that \propref{prop:small_clique}
holds. This implies that $|V(B)\setminus Z|\le (1+\eps)k$, hence
$|B| \le (1+2\eps)k$.

Next, we claim that there exists a pseudo-clique $C\in \mc C$, such that
$B$ contains $C$. To prove our claim, suppose, towards
contradiction, that $B$ does not contain any pseudo-clique from $\mc C$.
Then every vertex $u$ in $V(B)\setminus (Z\cup \{v_0\})$  must be outcast.
Indeed, if $u$ belongs to a pseudo-clique
$C'\in \mc C$, then $u$ is the smallest vertex in the unique block
that contains $C'$, hence $u\in Z_3\subseteq Z$, which is a
contradiction. Moreover, since every vertex in $V(B)\setminus (Z\cup
\{v_0\})$ is full in $B$, we have that $\Gamma$ is a graph with minimum
degree at least $(1-2\eps)k$ whose vertex set consists of only outcast
vertices. This fact, together with $|\Gamma| < (1+\eps)k$, implies that
the set $V(B)\setminus (Z\cup \{v_0\})$ forms a pseudo-clique in $G$
which is disjoint from all the other pseudo-cliques from $\mc C$.
But this contradicts the maximality of the union $\bigcup_{C\in \mc C}
C$, since we chose the collection of disjoint pseudo-cliques
that covers the maximum number of vertices possible. This contradiction
proves that $B$ contains exactly one pseudo-clique $C$ from $\mc C$.

It remains to show that $V(B)\subseteq Z \cup C$, or equivalently
$V(B)\setminus (C\cup Z) = \emptyset$. Suppose not. Clearly $v_0\in
Z_3\subseteq Z$ (because $v_0$ is the smallest vertex in $B$, and $B$
contains a pseudo-clique). Every vertex in $V(B)\setminus (Z\cup C)$
must be outcast and full in $B$, hence the graph $\Gamma':=
G'[C\cup (V(B)\setminus Z)]$ is a graph with minimum degree
$(1-4\eps)k$. By the assumption that \propref{prop:small_clique} holds,
we have that $|\Gamma'| < (1+ \frac{\eps}{2})k$. Hence the vertices
of $\Gamma'$ form a pseudo-clique in $G$, and if we replace $C$ by
$C':=C\cup (V(B)\setminus Z)$ (recall that $C'\setminus C$ consists only
of outcast vertices), we obtain a family of pseudo-cliques
whose union is larger than before, contradicting the maximality
of $\bigcup_{C\in \mc C} C$. This final contradiction establishes
the lemma.
\end{proof}

We turn to prove \lemref{lem:outcast_big}. One important fact that will be
used in subsequent double-counting arguments is the following consequence
of \obsref{obs:block_properties}: by removing the smallest vertex with
respect to $\le_T$ from each block in $\mc B$, we obtain a family of
pairwise vertex-disjoint graphs.

\begin{proof}[Proof of \lemref{lem:outcast_big}{}]
Assume, towards contradiction, that $\ell > 10^7\cdot\frac{n}{\eps k}$,
but $G_p$ does not \aas contain a cycle of length at least $k+1$.
The three lemmas \ref{lem:few_poor}, \ref{lem:many_poor}, and
\ref{lem:some_poor} combined imply that either $G_p$ \aas contains
a cycle of length at least $k+1$, or all the good blocks contain
a pseudo-clique inside, and the remaining vertices not in the pseudo-clique
are in $Z$. We can bound the number $t$ of vertices not in good blocks
as follows.

We claim that every block $B\in \mc B\setminus \mc B'$
of size $1 \le |B| < (1-4\eps)k$ must contain at least $\max\{1,|B| - 1\}$
vertices from $Z$ (and hence $B$ is necessarily not good). This is because
$B$ can not contain a pseudo-clique (its size is too small) and every
vertex of $V(B)\setminus Z$ which is not the smallest with respect to
$\le_T$ must be full in $B$ (see \corref{cor:not_z_full}).
However $B$ contains no full vertex, as $|B| < (1-4\eps)k$, therefore
$B$ contains at least $|B|-1$ vertices from $Z$. When $|B| = 1$, the
unique vertex in $B$ is isolated in $G'_p$, hence it belongs to $Z$
(and belongs to no other block in $\mc B$) thus proving our claim.

In particular, every non-good block $B$ contains at least $\min\lt\{|B|-1,
\frac{\eps}{10^6}|B|\rt\}$ vertices from $Z$. Furthermore, as it was
previously remarked, if we remove the smallest vertex from each block
in $\mc B$, we obtain a family of disjoint graphs. Hence, if $t_0$
denotes the number of blocks in $\mc B$ of size $1$, then
\begin{equation}
  \label{eqn:bound_not_good}
  |Z| \ge t_0 + \sum_{B\in \mc B}
  (|Z\cap B| - 1) \ge t_0 + \frac{\eps}{2\cdot 10^6} \sum_{\substack{%
  B \text{ not good} \\ |B| > 1}}
  |B|\ge \frac{\eps t}{2\cdot 10^6},
\end{equation}
thus $t \le \frac{2\cdot10^6}{\eps}\cdot |Z|$. By \eqref{eqn:z_small}
we obtain $t \le 10^6\cdot\lt(\frac{3n}{\eps k} + \frac{32}{\eps^2 pk}
\cdot \lt(\ell + \frac{n}{k}\rt)\rt)$.

Using \lemref{lem:some_poor}, we can estimate the number of
outcast vertices by adding the estimation of $|Z|$ in \eqref{eqn:z_small}
with our previous bound for $t$ in \eqref{eqn:bound_not_good} for $t$.
This is true because if an outcast vertex is in a good block, then it
must belong to $Z$. Hence we have $\ell \le t + |Z|$, which implies
that $\ell \le 10^7 \cdot \frac{n}{\eps k}$,
a contradiction that establishes \lemref{lem:outcast_big}.
\end{proof}

Next is \lemref{lem:outcast_small}.
\begin{proof}[Proof of \lemref{lem:outcast_small}{}]
Suppose that $\ell \le 10^7\cdot\frac{n}{\eps k}$, but $G_p$ does not
\aas have a cycle of length at least $k+1$. Clearly $|\mc C|$ is
roughly $\frac{n}{k}$, as the number of outcast vertices is $\ell
=o(n)$. Let $\mc B'$ be the sub-family of $\mc B$ consisting of the
good blocks. If we plug the inequality $\ell \le 10^7\cdot \frac{n}{\eps k}$
into the bound \eqref{eqn:z_small}, we obtain $|Z| \le \frac{1.1n}{k}$
since $\eps^2 p k \to \infty$ as $k\to \infty$. Moreover, using inequality
\eqref{eqn:bound_not_good}, we obtain that the number of vertices of $G'$
not in good blocks is at most $10^{7}\cdot \frac{n}{\eps k} = o(n)$.
Because of \lemref{lem:some_poor}, every member $B$ of $\mc B'$ contains
a pseudo-clique $C\in\mc C$ and $V(B)\setminus C \subseteq Z$, hence
$|\mc B'| \approx \frac{n}{k}$, or more specifically,
$0.99 n / k \le |\mc B'| \le 1.01n/k$. Furthermore, by corollaries
\ref{cor:tested_dfs} and \ref{cor:tested_block}, the total number of edges
in $Q_2\cup Q_3$ (recall that $Q_2$ is the set of edges tested by DFS, and
$Q_3$ is the set of edges tested by the block algorithm) is at most
$10^8 \cdot \frac{n}{\eps p k} = o(\eps n)$.

The number of blocks $B\in \mc B\setminus \mc B'$ having size
$|B| \ge (1-4\eps)k$ is $o(\frac{n}{k})$. This is because the total
number of vertices not in good blocks in $G'$ is at most $o(n)$.
Furthermore, every block $B\in \mc B$ of size $|B| < (1-4\eps)k$
has at least $\max\{1,|B|-1\}$ vertices in $Z$, as it was remarked
in the proof of \lemref{lem:outcast_big}. Thus the number of blocks in
$B\in \mc B \setminus \mc B'$ having size $|B| < (1-4\eps)k$
is at most $|Z|\le \frac{1.1n}{k}$. Combining these observations
together, we obtain $|\mc B| \le 3 |\mc B'|$.

The above implies the following four statements.

\begin{enumerate}[(i)] \smalllist
\item Fewer than $\frac{1}{5}\cdot|\mc B'|$ block in $\mc B'$ have
more than $5$ vertices from $Z$. This is a consequence of the inequality
$|Z| \ge \sum_{B\in \mc B} (|Z\cap V(B)|-1)$.

\item The number of blocks in $\mc B'$ having more than $5$
cut-vertices is less than $\frac{1}{5}\cdot |\mc B| \le \frac{3}{5}\cdot
|\mc B'|$. This is because the total number of cut-vertices is at most
$|\mc B| - 1$, and after the removal of the smallest vertex from each
block in $\mc B$, those blocks containing more than $5$ cut-vertices
still contribute at least $5$ to this total.

\item The vast majority of the blocks $B\in \mc B'$ are such that less
than $\eps k/3$ edges from $G'$ having one endpoint in $V(G')\setminus V(B)$
and the other being a non-cut-vertex of $B$. This is true since every such
edge is necessarily tested and belongs to $Q_2\cup Q_3$, and the total number
of edges in $Q_2\cup Q_3$ is $o(\eps n)$.

\item At least $0.99 |\mc B'|$ blocks in $\mc B'$ contain
pseudo-cliques that satisfy the condition \eqref{eqn:edgesoutside}
in \lemref{lem:edgesoutside}.

\end{enumerate}
Hence, there exists a block $B\in \mc B'$ satisfying the conditions
stated in (i)--(iv). Let $C\in \mc C$ be the unique
pseudo-clique contained in $B$.

It is time to incorporate the waste vertices back.
Let $N$ denote the union of $V(B)\cap Z$ with all cut-vertices
from $B$. By (i) and (ii), the set $N$ has size at most $10$ and clearly
$V(B)\subseteq C \cup N$. Let $F$ be the set of edges in $G$ connecting
$C\setminus N$ to a vertex outside of $C\cup N$. We prove that
$|F|\le\eps k$. In order to prove such inequality, we will apply
\lemref{lem:edgesoutside}. Let us recall the definitions of the sets
$D_1$, $D'_1$, $D_2$ and $\mc E$ in \eqref{eqn:edgesoutside}. The set $D_1$
consists of all the vertices in $C$ that have more than $\eps k$ neighbors
outside of $C$ in $G$. The subset $D'_1\subseteq D_1$ is the union of
$D_1\cap W$ with all vertices in $D_1\setminus W$ which lost more than
a $\frac{1}{100}$ proportion of its neighbors outside of $C$ after the
deletion of $W$. The set $D_2$ is the set of all vertices from $G$ not
in $C$ that have at least $\eps k$ neighbors in $C$. At last, $\mc E$ is
the set of all edges from $G$ connecting $C\setminus D_1$ to a vertex
outside not in $C\cup D_2$.

We claim that all vertices in $D_1\setminus D'_1$ are cut-vertices, and thus
$D_1\setminus D'_1 \subseteq N$. Every vertex in $D_1\setminus D'_1$
sends at least $0.99\eps k$ edges outside of $C$ in $G'$, in particular,
it also sends at least $0.99\eps k - |N| > \eps k / 3$ edges outside of $B$.
But by (iii), every vertex in $B$ that sends at least $\eps k / 3$
edges to the outside of $B$ must be a cut-vertex, hence
$D_1\setminus D'_1\subseteq N$.

Next, we claim that $D'_1=\emptyset$. By the discussion in the previous
paragraph, we have $D_1\setminus D'_1\subseteq N$, hence $|D_1\setminus D'_1|
\le 10$. The inequality \eqref{eqn:edgesoutside} states that $|D'_1|\le
|D_1|/100$, thus $|D_1\setminus D'_1| \ge 99 |D'_1|$. Hence we have
$|D'_1| < \frac{10}{99}<1$, which implies $D'_1 = \emptyset$.

Our next claim states that $D_2\setminus W \subseteq N$. To prove this,
let us first show that every vertex in $D_2\setminus W$ belongs to $B$.
A vertex in $D_2\setminus W$, sends at least $\eps k$ edges to $C$ in $G$,
and since $|C\cap W| < \eps k / 2$ (as otherwise $C$ would be completely
thrown away to the waste), this implies that every vertex in $D_2\setminus W$
has more than $\eps k / 2$ neighbors in $C\setminus W$ in the graph $G'$,
hence more than $\eps k / 2 - |N|$ neighbors in $C\setminus (W\cup N)$.
By (iii), any such vertex must belong to $B$, hence $D_2\setminus W
\subseteq V(B)$. On the other hand, since $D_2$ and $C$ are disjoint,
clearly we must have $D_2\setminus W \subseteq V(B)\cap Z \subseteq N$,
proving our claim.

Similarly to the proof of $D'_1=\emptyset$, let us now show that
$D_2\cap W = \emptyset$. Because $D_2\setminus W\subseteq N$, we have
$|D_2\setminus W| \le 10$, and by \eqref{eqn:edgesoutside}, we must
have $|D_2\setminus W| \ge 99 |D_2\cap W|$, therefore
$|D_2\cap W| < \frac{10}{99} < 1$, which implies $D_2\cap W = \emptyset$.

We turn to prove $|F| < \eps k$. Assume not. Because $D_1,D_2\subseteq N$,
we have $F\subseteq \mc E$. Moreover every vertex in $C\setminus D_1$ can
send at most $\eps k$ edges to the outside of $C$ in $G$, and every vertex
not in $D_2$ can send at most $\eps k$ edges to $C$ in $G$. Thus
$|\mc E\setminus F| \le |N\setminus D_1| \eps k + |N\setminus D_2|\eps k
\le 20\eps k$, hence $|F| \le |\mc E| \le |F| + 20 \eps k \le 21 |F|$.
By the inequality \eqref{eqn:edgesoutside}, we know that $|\mc E\setminus E(G')|\le
|\mc E|/100$, hence $|F\setminus E(G')|\le 21 |F|/100$. This implies that
there are at least $0.79\eps k>\eps k /3$ edges in $G'$ connecting a vertex
from $C\setminus N$ to a vertex outside of $C\cup N \supseteq V(B)$, which
contradicts (iii), hence $|F| < \eps k$. Therefore the pair $(C,N)$ satisfies
the statement of \lemref{lem:outcast_small}. This concludes the proof
of the lemma.
\end{proof}

\subsection{Step 6: finishing the proof}
\label{sec:step6}

It remains to prove \lemref{lem:good_pair}.

\begin{proof}[Proof of \lemref{lem:good_pair}{}]
If we remove the vertices from $N$ that have less than $\eps k$ neighbors
in $C\cup N$, we might increase the number of edges between $C\setminus N$
and $V(G)\setminus (C\cup N)$ to at most $\eps k + |N|\eps k \le 11\eps k$.
So we can assume that every vertex from $N$ has at least $\eps k$ neighbors
from $G$ in $C$, and that $e_G(C\setminus N, V(G)\setminus (C\cup N))\le
11\eps k$. Let $X = C \cup N$. From now on, we only deal with the graph $G[X]$.

Observe for the beginning that $|X|\ge k + 1$. Indeed, the minimum degree
in $G$ is at least $k$, and the number of edges between $C$ and
$V(G)\setminus X$ is less than the size of $C$. In the following,
we show that \aas $G_p[X]$ is Hamiltonian.

The general framework of the proof is to show some expansion properties
of $G_p[X]$ and then to deduce the Hamiltonicity from them.
There are several recent papers dedicated to or just using
Hamiltonicity of expanders, and the notion of expanders is slightly
different every time, depending on the setting it should be applied in.
Here we go with the notion used by Glebov and Krivelevich~\cite{GlKr13}:
a graph $H$ with the vertex set $[m]$ is called a $p'$-{\em expander},
if there exists a set $D\subset [m]$ such that $H$ and $D$ satisfy the
following properties:
\begin{itemize} \smalllist
\item $|D|\le m^{0.09}.$
\item The graph $H$ does not contain a non-empty path of length at most
$\frac{2\log m}{3\log\log m}$ such that both of its (possibly identical)
endpoints lie in $D$.
\item For every set $S\subset [m]\setminus D$ of size $|S|\le \frac{1}{p'}$,
its neighborhood satisfies $|N(S)|\ge \frac{m p'}{1000}|S|$.
\end{itemize}

Let us denote $F=G[X]$, and let $m=|X|$ be its order. Furthermore,
let us define for convenience $p'=\frac{\log m + \log\log m}{m}$.
We first show that for every $f_1, f_2=\omega_m(1),\, f_1<f_2<\log\log m$,
every graph $H$ satisfying $F_{p_1}\subseteq H\subseteq F_{p_2}$ with
$p_i=\frac{\log m + \log\log m +f_i}{m}$ is \aas a $p'$-expander.
(Notice that we are coupling $F_{p_1}$ and $F_{p_2}$, so that $F_{p_1}
\subseteq F_{p_2}$)
Indeed, let us fix $D=\{v\in X:\, d_{F_{p_1}}(v)<m p'/100\}$ to be the
set of all vertices from $X$ with degree less than $m p'/100$ in $F_{p_1}$.
The proof of the first property is similar to the proof of Claim~4.3
in~\cite{BSKrSu11},  and the second property is shown to hold similarly
to Claim~4.4 in~\cite{BSKrSu11}. Finally, the proof of the third bullet
follows the lines of the corresponding proof in Lemma~10 in~\cite{GlKr13}.
Furthermore, observe that \lemref{lem:chernoff} guarantees us that \aas
every vertex from $C$ has degree at least two in $H$, and for the vertices
in $N$ this also holds \aas by \lemref{lem:chernoff}. Hence, the random
graph $H$ is \aas a $p'$-expander with minimum degree at least $2$.
Applying Lemma~11 from~\cite{GlKr13}, we see that $H$ is either Hamiltonian
or has quadratically many boosters.

With this statement in our toolbox, the proof is similar to the proof of
\propref{prop:good_clique}. We fix $p_1$ such that $p-p_1=\omega_m(1)$,
and let $p_2=p$. We start with $H=F_{p_1}$ and successively add random
edges to its edge set until we obtain $F_{p}$. We update the set of
boosters after each new edge. Every such edge has at least constant
probability to be a booster for the current $H$ as long as $H$ is not
Hamiltonian. Every added edge that is a booster increases the length of
the longest cycle in the current graph by at least one, or makes it
Hamiltonian. Therefore, after at most $k$ added boosters, the process
would end with a Hamiltonian graph.
On the other hand, the total number of added edges is a binomial random
variable $|E(F_p)|-\lt|E\lt(F_{p_1}\rt)\rt|$ with $\binom{m}{2}$ trials and
probability $p-p_1=\omega_m(1)$. By \lemref{lem:chernoff}, with probability
at least $1-\exp(-m)$, the number of new edges that are added to obtain
$F_{p}$ from $F_{p_1}$ is $\omega_m(m)$.
Hence, \lemref{lem:chernoff} guarantees us that \aas we get sufficiently
many boosters to make the graph $F_{p}$ Hamiltonian, proving the lemma.
\end{proof}

\section{Concluding remarks and open questions}
\label{sec:conclusion}
In this paper, we studied random subgraphs of graphs with large minimum degree. Our goal was to
extend classical results on random graphs to a more general model, where we replace the host graph
by a graph with large minimum degree. We determined the threshold probability for having cycle of length at least
$k+1$ in the random subgraph of graph with minimum degree at least $k$, showing that the assertion about
Hamiltonicity of $\gnp(k+1,p)$ can be extended to this setting.

We believe that there are further interesting statements that one can deduce from our proof.
One of them is the bipartite version of \thmref{thm:main}. Namely, that in a
bipartite graph with minimum degree at least $k$,
the random subgraph (with the same probability as in this paper)
\aas contains a cycle of length at least $2k$. However, since the paper
is already quite long, we do not check all the technical details needed for
the proof of this statement.

Another fact that can be shown similarly to \thmref{thm:main}
is as follows. Let $G$ be a graph with minimum degree at least $k$,
and fix a constant $c$. If $p=p(k)\ge\frac{\log k + \log \log k + c}{k}$,
then $G_p$ contains a cycle of length at least $k+1$ with probability at
least $e^{-e^{-c}}-o(1)$. This particular statement is an analog of the
well-known result on the probability of $\gnp(k+1,p)$ being Hamiltonian
in the range of $p$ where the probability of having one vertex of degree
at most one is a constant (see, e.g., \cite{BolBook}). The only difference
in the proof compared to \thmref{thm:main} would be the proof of the
corresponding version of \lemref{lem:good_pair}, since this is the only
place where we use the additional summand $\omega(1)$ in the definition of
$p$.

One natural question is to determine whether the results of this paper, as well as several previous
ones on this topic, hold if one weakens the condition of minimum degree of the
host graph. One possibility here would be to only require the host graph
$G$ to have \emph{average degree} at least $k$.
Does this still guarantee cycles of length $(1-o(1))k$ and $k+1$ in $G_p$,
for the same value of $p$ as in~\cite{KrivLeeSud} and in this paper?

Finally, it would be interesting to find more monotone properties $\mc P$ for
which the threshold probability in the binomial random graph model
is the smallest among all host graphs of given minimum degree.
Formally, these are the properties $\mc P$ such that if $\gnp(n,p)$ \aas satisfies $\mc P$,
then this holds \aas  also for a random subgraph $G_p$ of a graph $G$ with minimum degree
at least $n-1$.

\end{document}